\documentclass[12pt,reqno]{amsart}
\usepackage{amsmath,amsfonts,amsthm,amsopn,amssymb}
\usepackage{cite,marginnote}
\usepackage{color,enumitem,graphicx}
\usepackage[colorlinks=true,urlcolor=blue,
citecolor=red,linkcolor=blue,linktocpage,pdfpagelabels,
bookmarksnumbered,bookmarksopen]{hyperref}
\usepackage[english]{babel}

\usepackage[left=2.9cm,right=2.9cm,top=2.8cm,bottom=2.8cm]{geometry}
\usepackage[hyperpageref]{backref}

\numberwithin{equation}{section}

\makeindex
\newtheorem{theorem}{Theorem}[section]
\newtheorem{lemma}{Lemma}[section]
\newtheorem{corollary}{Corollary}[section]
\newtheorem{proposition}{Proposition}[section]
\newtheorem{remark}{Remark}[section]
\newtheorem{example}{Example}[section]
\newtheorem{definition}{Definition}[section]

\newcommand{\ud}{\mathrm{d}}

\newcommand{\iy}{\infty}

\newcommand{\s}{\section}

\newcommand{\la}{\lambda}

\newcommand{\R}{\mathbb R}

\newcommand{\C}{\mathbb C}

\newcommand{\rg}{\rightarrow}

\newcommand{\lan}{\langle}
\newcommand{\ran}{\rangle}

\newcommand{\vp}{\varphi}

\newcommand{\bt}{\begin{theorem}}
\newcommand{\et}{\end{theorem}}
\newcommand{\bl}{\begin{lemma}}
\newcommand{\el}{\end{lemma}}
\newcommand{\bd}{\begin{definition}}
\newcommand{\ed}{\end{definition}}
\newcommand{\bc}{\begin{corollary}}
\newcommand{\ec}{\end{corollary}}
\newcommand{\bp}{\begin{proof}}
\newcommand{\ep}{\end{proof}}
\newcommand{\bx}{\begin{example}}
\newcommand{\ex}{\end{example}}
\newcommand{\bi}{\begin{exercise}}
\newcommand{\ei}{\end{exercise}}
\newcommand{\bo}{\begin{proposition}}
\newcommand{\eo}{\end{proposition}}
\newcommand{\br}{\begin{remark}}
\newcommand{\er}{\end{remark}}
\newcommand{\be}{\begin{equation}}
\newcommand{\ee}{\end{equation}}
\newcommand{\ba}{\begin{align}}
\newcommand{\ea}{\end{align}}
\newcommand{\bn}{\begin{enumerate}}
\newcommand{\en}{\end{enumerate}}
\newcommand{\bg}{\begin{align*}}
\newcommand{\bcs}{\begin{cases}}
\newcommand{\ecs}{\end{cases}}

\newcommand{\bean}{\begin{eqnarray*}}
\newcommand{\eean}{\end{eqnarray*}}

\def\Proof{\noindent{\bf Proof}\quad}

\title[Existence and multiplicity of sign-changing standing waves]
{Existence and multiplicity of sign-changing standing waves for a gauged nonlinear Schr\"{o}dinger equation
in $\R^2$}

\author[Z.S.~Liu]{Zhisu Liu}
\author[Z.G.~Ouyang]{Zigen Ouyang}
\author[J.J.~Zhang]{Jianjun Zhang}

\address[Z.S.~Liu]{\newline\indent School of Mathematics and Physics,
\newline\indent University of South China,
\newline\indent 421001, Hengyang, Hunan, P. R. China}

\email{\href{mailto:liuzhisu183@sina.com}{liuzhisu183@sina.com}}

\address[Z.G.~Ouyang]{\newline\indent School of Mathematics and Physics,
\newline\indent University of South China,
\newline\indent 421001, Hengyang, Hunan, P. R. China
}
\email{\href{mailto:zigenouyang@163.com}{zigenouyang@163.com}}

\address[J.J.~Zhang]{\newline\indent College of Mathematics and Statistics
\newline\indent
Chongqing Jiaotong University
\newline\indent
Chongqing 400074, PR China
\newline\indent and
\newline\indent Dip. di Scienza e Alta Tecnologia
\newline\indent
Universit\`{a} degli Studi dell'Insubria
\newline\indent
via G.B. Vico 46, 21100 Varese, Italy}
\email{\href{mailto:jianjun.zhang@uninsubria.it}{jianjun.zhang@uninsubria.it}}

\thanks{Z.S.~Liu was partially supported by the NSFC (Grant No.
11271115) and NSF of Hunan Province (No. 2017JJ3265), the Science and Technology Plan Project of Hengyang
 City(2017KJ183), Research Foundation of Education Bureau of Hunan Province(16A031). The corresponding author is J.J.~Zhang, who was partially supported by NSFC(No. 11871123).}


\begin{document}

\begin{abstract}
We are concerned with sign-changing solutions of the following gauged nonlinear Schr\"{o}dinger
equation in dimension two including the so-called Chern-Simons term
\begin{align*}
\left\{
  \begin{array}{ll}
    -\triangle {u}+\omega u+\left(\frac{h^2(|x|)}{|x|^2}+\int_{|x|}^{+\infty}\frac{h(s)}{s}u^2(s){\rm ds}\right) u =\lambda|u|^{p-2}u& \mbox{in}\,\,\R^2, \\
    u(x)=u(|x|)\,  \in\, H^1(\R^2),
  \end{array}
\right.
\end{align*}
where $\omega,\lambda>0$, $p\in(4,6)$ and
$$
h(s)=\frac{1}{2}\int_0^s\tau u^2(\tau)d\tau.
$$
Via a novel perturbation approach and the method of invariant sets of descending flow,
we investigate the existence and multiplicity of sign-changing solutions.
Moreover, {\it energy doubling} is established, i.e., the energy of sign-changing solution $w_\lambda$
is strictly larger than twice that of the ground state energy for $\lambda>0$ large.
Finally, for any sequence $\lambda_n\rightarrow\infty$ as $n\rightarrow\infty$,
up to a subsequence, $\lambda_n^{\frac{1}{p-2}}w_{\lambda_n}\rg w$
strongly in $H_{rad}^1(\R^2)$ as $n\rightarrow\infty$, where $w$ is a sign-changing solution of
$$
-\triangle {u}+\omega u=|u|^{p-2}u,\,\,u\in H_{rad}^1(\R^2).
$$
\end{abstract}

\maketitle

\s{Introduction and Main Results}
\renewcommand{\theequation}{1.\arabic{equation}}
\subsection{Gauged Schr\"{o}dinger equations}
Consider the following planar gauged
nonlinear Schr\"{o}dinger system
\begin{equation}\label{eqn:cs1}
\left\{
  \begin{array}{ll}
    -iD_0\phi+ (D_1D_1+D_2D_2)\phi =-\lambda|\phi|^{p-2}\phi, \\
    \partial_0A_1-\partial_1A_0=-\text{Im}(\bar{\phi}D_2\phi),\\
    \partial_0A_2-\partial_2A_0=\text{Im}(\bar{\phi}D_1\phi),\\
   \partial_1A_2-\partial_2A_1=-\frac{1}{2}|\phi|^2,
  \end{array}
\right.
\end{equation}
where $i$ denotes the imaginary unit, $\partial_0=\frac{\partial}{\partial t}$,
$\partial_1=\frac{\partial}{\partial x_1}$, $\partial_2=\frac{\partial}{\partial x_2}$ for
$(t,x)\in\R^{1+2}$, $x=(x_1,x_2)$, $\phi:\,\R^{1+2}\rightarrow\C$ is the complex scalar field,
$A_j:\,\R^{1+2}\rightarrow\R$ is the gauge field,
$D_j=\partial_j+iA_j$ is the covariant derivative
for $j=0,1,2$, and $\lambda$ is a positive constant representing the strength of interaction potential.
One feature of this model is that this system was proposed in  consists of Schr\"{o}dinger equations
augmented by the gauge field. And it is particular important in studying the high-temperature
superconductor, fractional quantum Hall effect and Aharovnov-Bohm scattering. For more details about system (\ref{eqn:cs1}), we refer the readers to \cite{Dunne95,Hagen84,Hagen85}.
Moreover, system (\ref{eqn:cs1}) is invariant under the following gauge transformation
\begin{equation}\label{eqn:invariant}
\phi\rightarrow\phi e^{i\chi}, \quad A_j\rightarrow A_j-\partial_j\chi,
\end{equation}
 where $\chi:\,\R^{1+2}\rightarrow\R$ is an arbitrary $C^\infty$ function. This system
 was first studied in \cite{Jackiw1,Jackiw2,Jackiw92}. The initial value problem of such a system
as well as global existence and blow-up,
has also been addressed in \cite{Berge95,Huh09,Huh12,Liu14} for
the case $p=4$. we can also see \cite{Liu16} for a global existence
result in the defocusing case, and \cite{Chen14} for a uniqueness result to the infinite radial
hierarchy.

The existence of stationary states to system (\ref{eqn:cs1}) with general $p>2$ has been studied
recently in \cite{Byeon12}. By using the ansatz
$$
\aligned
&\phi(t,x)=u(|x|)e^{i\omega t},\quad A_0(t,x)=A_0(|x|),\\
&A_1(t,x)=\frac{x_2}{|x|^2}h(|x|),\quad A_2(t,x)=-\frac{x_1}{|x|^2}h(|x|),
\endaligned
$$
Byeon, Huh and Seok found in \cite{Byeon12} that $u$ solves the equation
\begin{equation}\label{eqn:cs2}
-\triangle u+(\omega+\xi)u+\left(\frac{h^2(|x|)}{|x|^2}+\int_{|x|}^{+\infty}\frac{h(s)}{s}u^2(s){\rm d}s\right)u
=\lambda|u|^{p-2}u,
\end{equation}
where $h(s)=\frac{1}{2}\int_0^s\tau u^2(\tau){\rm d}\tau$ and $\xi$ is an integration constant of $A_0$, which takes the form
$$
A_0(r)=\xi+\int_{r}^{+\infty}\frac{h(s)}{s}u^2(s){\rm d}s.
$$
As mentioned in \cite{Byeon12},  taking $\chi =ct$ in the gauge invariance (\ref{eqn:invariant}), we
obtain another stationary solution for any given stationary solution;
the functions $u(x), A_1(x), A_2(x)$ are preserved,
and
$$
\omega\rightarrow\omega+c,\quad A_0(x)\rightarrow A_0(x)-c.
$$
That is to say, the constant $\omega+\xi$ is a gauge invariant of the stationary solutions of
the problem.  Therefore, we can take $\xi=0$ in what follows, that is,
$$
\lim\limits_{|x|\rightarrow\infty}A_0(x)=0,
$$
which was indeed assumed in \cite{Jackiw92}. In this situation, equation (\ref{eqn:cs2}) becomes
\begin{equation}\label{eqn:cs3}
-\triangle u+\omega u+\left(\frac{h^2(|x|)}{|x|^2}+\int_{|x|}^{+\infty}\frac{h(s)}{s}u^2(s){\rm d}s\right)u
=\lambda|u|^{p-2}u,\quad x\in\R^2.
\end{equation}
It is shown in \cite{Byeon12} that (\ref{eqn:cs3}) is indeed variational and the associated energy functional
$I_\lambda:\,H_{rad}^1(\R^2)\rightarrow\R$ is given by
$$
I_\lambda(u)=\frac{1}{2}\int_{\R^2}(|\nabla u|^2+\omega u^2){\rm d}x+\frac{1}{2}\int_{\R^2}\frac{u^2}{|x|^2}
h^2(|x|){\rm d}x-\frac{\lambda}{p}\int_{\R^2}|u|^{p}{\rm d}x.
$$
Here $H_{rad}^1(\R^2)$ denotes the subspace of radially symmetric functions in $H^1(\R^2)$ with
the inner product and norm
$$
\langle u,v\rangle=\int_{\R^2}(\nabla u\nabla v+\omega uv){\rm d}x,
\quad \|u\|:=\left(\int_{\R^2}(|\nabla u|^2+\omega u^2){\rm d}x\right)^{\frac{1}{2}}.
$$
For simplicity, in what follows, denote
$$
B(u):=\frac{1}{2}\int_{\R^2}\frac{u^2}{|x|^2}h^2(|x|){\rm d}x,
$$
then
$B\in C^1(H_{rad}^1(\R^2),\R)$ and
$$
\langle B'(u),\varphi\rangle=\int_{\R^2}\left(\frac{h^2(|x|)}{|x|^2}+\int_{|x|}^{+\infty}\frac{h(s)}{s}u^2(s)\,\ud s\right)u(x)\varphi(x)\,\ud x,\,\,\forall\,\,\varphi\in H_{rad}^1(\R^2).
$$

For any $u,\varphi\in H_{rad}^1(\R^2)$, we have
\begin{equation}\label{eqn:function}
\aligned
\langle I'_\lambda(u),\varphi\rangle=&\int_{\R^2}(\nabla u\nabla\varphi+\omega u\varphi){\rm d}x
-\lambda\int_{\R^2}|u|^{p-2}u\varphi {\rm d}x+\langle B'(u),\varphi\rangle
\endaligned
\end{equation}
And any critical point of $I_\lambda$ in $H_{rad}^1(\R^2)$ is a weak solution
of equation (\ref{eqn:cs3})(see \cite{Byeon12}).

\subsection{Motivation}

In the following, we summarize some relative results for equation (\ref{eqn:cs3}).
In contrast with Schr\"odinger equations, equation (\ref{eqn:cs3}) is {\it nonlocal},
that is, it not a pointwise identity
with the appearance of the Chern-Simons term
$$
\left(\frac{h^2(|x|)}{|x|^2}+\int_{|x|}^{+\infty}\frac{h(s)}{s}u^2(s){\rm d}s\right)u.
$$
This nonlocal term causes some mathematical difficulties and makes the problem rough and particularly interesting. One of main obstacles is
the boundedness of Palais-Smale sequences if one tries to use directly the Mountain Pass Theorem
to get critical points of $I_\lambda(u)$ in $H_{rad}^1(\R^2)$. For $p\geq6$,
it is standard to show that the Palais-Smale condition holds for $I_\lambda(u)$ in $H_{rad}^1(\R^2)$,
see \cite{Huh12}.
For $p>4$, the Euler-Lagrange functional is unbounded
from below and exhibits a mountain-pass geometry.
However, it seems hard to prove that the Palais-Smale condition holds for $p \in(4,6)$.
 This problem was bypassed in \cite{Byeon12} by using a constrained
minimization taking into account both the Nehari and Pohozaev identities.
Moreover, in \cite{Byeon12}, static solutions can be
found for the special case $p=4$ by passing to a self-dual equation, which leads to
a Liouville equation that can be solved explicitly.
But those solutions are positive. In the same paper, the case $p\in(2, 4)$ also
was considered by using a minimization argument on a $L^2$-sphere.
By analysing the corresponding limit equation and investigating
the geometry of the Euler-Lagrange functional, Pomponio and Ruiz \cite{Pomponio15}
 proved the existence and nonexistence of positive solutions to equation (\ref{eqn:cs3}) for different
value of $\lambda$ when $p\in(2, 4)$. They showed that there exists
a threshold value $\omega_0>0$ such that the corresponding energy functional
is bounded from below if $\omega> \omega_0$, and not for $\omega\in(0,\omega_0)$.
Moreover, they gave an explicit expression of $\omega_0$ as follows
$$
\omega_0=\frac{4-p}{2+p}3^\frac{p-2}{2(4-p)}2^{\frac{2}{4-p}}\left(\frac{m^2(4+p)}{p-2}\right)^{-\frac{p-2}{2(4-p)}},
$$
where

$$
m=\int_{-\infty}^\infty\left[\frac{2}{p}\coth^{2}\left(\frac{p-2}{2}r\right)\right]^{\frac{2}{2-p}}{\rm d}r
$$
is the unique positive even solution of the problem
$$
-v''+v=v^{p-1}\quad\text{in}\,\,\R.
$$
For $p\in(2, 4)$, boundary concentration of solutions to (\ref{eqn:cs3}) was considered by Pomponio and Ruiz \cite{Pomponio15-} in the case of bounded domains. As for more general nonlinearities, there have been some works in which the authors tried to weaken or eliminate the Ambrosetti-Rabinowitz condition. In the sprit of the Berestycki-Lions conditions in \cite{Berestycki83,Berestycki83-}, Cunha, d'Avenia, Pomponio and Siciliano \cite{Cunha15} explored the Chern-Simons-Schr\"{o}dinger system with a nonlinear term of Berestycki,
Gallou\"{e}t and Kavian type in \cite{Berestycki84}. In this aspect, we also would like to cite \cite{Wang14}.
On the other hand, Li and Luo \cite{Li17} also investigated the nonlocal equation (\ref{eqn:cs3}) in the mass-critical
case: $p=4$ and mass-supercritical case: $p>4$, for instance, the existence, $H^1(\R^2)$-bifurcation and multiplicity of
normalized solutions. The existence of stationary states
with a vortex point has also been considered in \cite{Byeon16,Jiang16}.
For more existence results on the nonlinear Chern-Simons-Schr\"{o}dinger equation (\ref{eqn:cs3}),
we refer the readers to \cite{Yuan15,Luo18,Seok15,Ji17} and the references therein.
\subsection{Main results}
To the best of our knowledge, in the literature, there are just few result on the existence
of sign-changing solutions for Chern-Simons-Schr\"{o}dinger systems.
Combining constraint minimization method and quantitative deformation lemma,
Li, Luo and Shuai \cite{Li17-} proved the existence and asymptotic behavior of the
least energy sign-changing solutions to (\ref{eqn:cs3}) when $p>6$.
Moreover, {\it energy doubling} for the sign-changing solutions was obtained. Meanwhile, Deng, Peng and Shuai \cite{Deng17} also established the existence and
asymptotic behavior of nodal solutions to (\ref{eqn:cs3}) when $p>6$.
Precisely, they obtained the existence of a sign-changing solution, which changes
signs exactly $k$ times for any $k\in \mathbb{N}$. Their procedure of arguments is to transform
the original problem to solving a system of $(k+1)$ equations with $(k+1)$ unknown functions
$u_i$ with disjoint supports. Then the nodal solution is constructed through gluing $u_i$ by
matching the normal derivative at each junction point. We highlight
that $p>6$ plays a crucial role in \cite{Li17-,Deng17}.
So, one question still remains open in these results above:
\begin{center}
{\it Does problem (\ref{eqn:cs3}) admit sign-changing solutions in the case $ p\in(4,6)$ ?}
\end{center}
The main interest of the present paper is to give an affirmative answer to this open question.
\vskip0.1in

\begin{theorem}{\bf(Existence)}\label{Thm:existence}
If $p\in(4,6)$, then for any $\lambda>0$, equation (\ref{eqn:cs3})
admits at least one least energy sign-changing solution in $H_{rad}^1(\R^2)$.
\end{theorem}

Since the Chern-Simons term
$$
\left(\frac{h^2(|x|)}{|x|^2}+\int_{|x|}^{+\infty}\frac{h(s)}{s}u^2(s){\rm d}s\right )u
$$
is involved, equation (\ref{eqn:cs3}) is quite different from the following
scalar field equation
\begin{equation}\label{eqn:ellptic}
-\triangle {u}+\omega u=|u|^{p-2}u,
\end{equation}
which does not depend on the nonlocal term any more. In the literature,
problem (\ref{eqn:ellptic}) has attracted a considerable attention since 1970s.
In \cite{Berestycki83}, the so-called Berestycki and Lions was introduced to
guarantee the existence of ground state solutions to a general scalar
field equation $-\triangle u =f(u)$. Moreover, the Berestycki and Lions
condition is almost necessary. At this aspect, we also refer to
\cite{Berestycki83-,Berestycki84}. With some suitable condition on $f$,
ground state solutions have a fixed sign in general. Another topic of particular
interest is sign-changing solutions, which can be sought by several variational approaches,
for instance, Nehari manifold method \cite{Bartsch01,Cao88,Cerami86},
Morse theory \cite{Chang04}, heat flow method \cite{Chang03} and so on.
Over the years, combining minimax methods, the method of invariant sets of descending flow has
been a powerful tool in finding sign-changing
solutions of elliptic problems. For more progress in this aspect,
we refer to \cite{Bartsch04,Bartsch05,Liusun,Liuwang,Zou} and the references therein.

\begin{remark} In contrast to problem (\ref{eqn:ellptic}),
the Chern-Simons term makes equation (\ref{eqn:cs3}) tough. In the following,
we summarize some difficulties caused by the non-locality in seeking
sign-changing solutions.

\begin{itemize}
\item[(1)] In finding sign-changing solutions of (\ref{eqn:ellptic}),
a crucial ingredient is the following decomposition: for any $u\in H^1(\R^2)$,
\begin{equation}\label{eqn:fenjie1}
I(u)=I(u^+)+I(u^-),
\end{equation}
\begin{equation}\label{eqn:fenjie2}
\langle I'(u), u\rangle=\langle I'(u^+), u^+\rangle,\quad\langle I'(u), u^-\rangle=\langle I'(u^-), u^-\rangle
\end{equation}
where $I$ is the energy functional associated to (\ref{eqn:ellptic}) and defined by
$$
I(u):=\frac{1}{2}\int_{\R^2}(|\nabla u|^2+\omega u^2){\rm d}x
-\frac{1}{p}\int_{\R^2}|u|^p{\rm d}x.
$$
However, in applying the method of invariant sets of descending flow to equation (\ref{eqn:cs3}), due to the Chern-Simons term, one of the main obstacles is that the following decompositions
\be\label{eqn:fenjie3}
I_\lambda(u)=I_\lambda(u^+)+I_\lambda(u^-),\,\,\,\langle I_\lambda'(u), u^\pm\rangle=\langle I_\lambda'(u^\pm), u^\pm\rangle
\ee
do not hold any more in general.
Motivated by \cite{Liuz16}, via a sign-changing critical point theorem due to J. Liu, X. Liu and Z.-Q. Wang \cite{Liuj15}, we attempt to seek sign-changing solutions for equation (\ref{eqn:cs3}). Noting that,
the spitting decomposition of the Chern-Simons term is more complicated than that of the nonlocal term
involved in Schr\"odinger-Poisson systems(see \cite{Liuz16}), we need some new trick.
Thus,
the variational framework in \cite{Liuj15,Liuz16} is not directly applicable due to
changes of the geometric nature of the energy functional $I_\lambda$.
\item[(2)]  The effect of the nonlocal term
$\left(\frac{h^2(|x|)}{|x|^2}+\int_{|x|}^{+\infty}\frac{h(s)}{s}u^2(s){\rm d}s\right)u$ results in two difficulties.
First, it seems much complicated to find a similar auxiliary operator A(see \cite[Section 4]{Liuz16}),
 which plays a crucial role in constructing invariants sets of a descending flow associated
 with equation (\ref{eqn:cs3}). A similar difficulty also arises in
 seeking sign-changing solutions of
 the Schr\"odinger-Poisson systems
\begin{equation}\label{SP} \left\{
\begin{array}{ll}
-\Delta u+V(x)u+\phi u=|u|^{p-2}u&\mbox{in}\ \R^3,\\
-\Delta\phi=u^2&\mbox{in}\ \R^3,
\end{array}
\right.
\end{equation}
where $p\in(3,4)$. In \cite{Liuz16}, the authors overcome this difficulty
for $p\in(3,4)$ by adding a higher order local nonlinear term than $4$ and a coercive
condition about $V$. Second, the fact that $p\in(4,6)$ makes
tough to get the boundedness of (PS) sequences. In \cite{Liuz16}, the authors
recovered such boundedness due to the coercivity condition of $V$. However, without the coercivity condition,
the method in \cite{Liuz16} will be not applicable any more. In this paper,
we develop particularly \emph{a perturbation approach} by adding another nonlocal perturbation. Then we can obtain sign-changing solutions of the perturbed problems. By passing to the limit, sign-changing solutions of the original equation
(\ref{eqn:cs3}) are obtained. We should point out that, this perturbation approach also works for system (\ref{SP}) when $p\in(3,4)$ without any coercive conditions. Moveover, a similar perturbation approach also can be found in \cite{Liu18}, but does not work in the present paper. In fact, an additional perturbation term is also involved in this paper(see Section 3.1). In some sense,  it can  be seen as a generalization of that in \cite{Liu18}.
\end{itemize}
\end{remark}

Our another purpose is to establish {\it energy doubling} of sign-changing solutions
to equation (\ref{eqn:cs3}) with $p\in(4,6)$. We use the {\it local} equation (\ref{eqn:ellptic}) to illustrate {\it energy doubling}. Set
$$
\mathcal{N}:=\{u\in H_{rad}^1(\R^2)\setminus\{0\}:\,\,\langle I'(u),u\rangle=0\},
$$
and
\begin{equation} \label{eqn:p-ground}
c_0=\inf\{I(u):\,u\in\mathcal{N}\}.
\end{equation}
Obviously, for any sign-changing solution $w\in H_{rad}^1(\R^2)$ of equation (\ref{eqn:ellptic}), one can get that
\begin{equation} \label{eqn:p-ground-}
I(w)=I(w^+)+I(w^-)\geq 2c_0.
\end{equation}
In \cite{Weth11}, Ackermann and Weth proved that for each sign-changing solution
$w$ of (\ref{eqn:ellptic}), it holds true that
$$
I(w)>2c_0,
$$
which is called in the literature that $w$ satisfies
\lq\lq energy doubling\rq\rq(see also \cite{Weth06}).
Now, we give the analogue of {\it energy doubling} for problem (\ref{eqn:ellptic}) as follows.

\begin{definition}
Let $w_\lambda\in H_{rad}^1(\R^2)$ be a sign-changing solution of equation (\ref{eqn:cs3}),
$w_\lambda$ satisfies energy doubling, if $I_\lambda(w_\lambda)>2c_\lambda$, where
\begin{equation} \label{eqn:p-ground1}
c_\lambda=\inf\{I_\lambda(u),u\in\mathcal{N}_\lambda\},
\quad\mathcal{N}_\lambda:=\{u\in H_{rad}^1(\R^2)\setminus\{0\}:\,I_\lambda'(u)=0\}.
\end{equation}
\end{definition}
Let $w_\lambda\in H_{rad}^1(\R^2)$ be a sign-changing solution
of equation (\ref{eqn:cs3}).
Since the interaction of the positive and negative parts of solutions can not be neglected,
it is known that
\begin{equation} \label{eqn:p-ground+}
w_\lambda^{\pm}\not\in \mathcal{N}_\lambda.
\end{equation}
Thus, a natural open question is that whether energy doubling holds or not.
Generally speaking, it is even not easy to compare $I_\lambda(w_\lambda)$ with $c_\lambda$.
By using an approximation procedure, we give a partial answer for such an open problem,
that is, energy doubling holds if $\lambda>0$ large.
Precisely, we have the following result.
\begin{theorem}{\bf(Energy doubling)}\label{Thm:ground}
If $p\in(4,6)$, then there exists $\lambda^*>0$ such that, for $\lambda>\lambda^*$, any sign-changing solutions
$w_\lambda$ of problem (\ref{eqn:cs3})
in $H_{rad}^1(\R^2)$ satisfies energy doubling.
Furthermore, for any sequence $\{\lambda_n\}$ with $\lambda_n\rightarrow+\infty$
as $n\rightarrow\infty$, up to a subsequence,
$\lambda_n^{\frac{1}{p-2}}w_{\lambda_n}\rightarrow w$ in $H_{rad}^1(\R^2)$,
where $w\in H_{rad}^1(\R^2)$ is a sign-changing solution of (\ref{eqn:ellptic}).
\end{theorem}

Our last purpose in the present paper is to look for infinitely many
sign-changing solutions of equation (\ref{eqn:cs3}) when $\lambda$ is small enough.
\begin{theorem}{\bf(Multiplicity)}\label{Thm:many}
If $p\in(4,6)$, then there exists $\lambda_*>0$ such that, for any $\lambda\in(0,\lambda_*)$
 equaiton (\ref{eqn:cs3}) has infinitely many sign-changing solutions in $H_{rad}^1(\R^2)$.
\end{theorem}
For the proof of Theorem \ref{Thm:many}, we adopt the perturbation approach
established in \cite{Liuj15,Liuz16}. Based on a new inequality of Sobolev type developed
by Byeon, Huh and Seok in \cite{Byeon12} and a perturbed term growing faster than $4$,
the boundedness of (PS)
sequences can be established. Via the method of invariants sets of a descending
flow and the limiting argument,
we establish the existence of multiple sign-changing solutions.


This paper is organized as follows. In Section 2, some
notations and preliminaries are given. Section 3 is devoted to proving
Theorem \ref{Thm:existence} via a perturbation argument and the method of invariant
sets of decreasing flow. Section 4 is denoted to proving Theorem \ref{Thm:ground}.
Moreover, the asymptotic behavior is established as the parameter $\la$ goes to $+\iy$.
Finally, Section 5 is denoted to proving Theorem \ref{Thm:many}.


\s{Preliminaries and Functional Setting}
\renewcommand{\theequation}{2.\arabic{equation}}
 Let us fix some notation and give some preliminary results. For every $1\leq s\leq +\infty$, we denote by $\|\cdot\|_s$ the usual
norm of the Lebesgue space $L^s(\R^2)$. $c,C$ denote (possibly different) positive constants which may change from line to line.

In the following, we introduce some properties of $A_i(i=0,1,2)$, $B$ and $h$,
which will be used later frequently.
\begin{lemma}[\cite{Byeon12}]\label{Lem:feijubu}
If $u\in H_{rad}^1(\R^2)\cap L_{loc}^\infty(\R^2)$, then $A_0\in L^\infty(\R^2)$. Furthermore,
if $u\in H_{rad}^1(\R^2)\cap C(\R^2)$, then $A_0,A_1$ and $A_2$ are in $H_{rad}^1(\R^2)\cap C^1(\R^2)$.
\end{lemma}
\begin{lemma}[\cite{Byeon12}]\label{Lem:bu}
\begin{itemize}
\item[\rm (1)] If $u_n\rightharpoonup u$ weakly
in $H_{rad}^1(\R^2)$ as $n\rightarrow+\infty$, then
\begin{itemize}
\item[\rm (i)] $\lim\limits_{n\rightarrow+\infty}B(u_n)=B(u)$;
\item[\rm (ii)]  $\lim\limits_{n\rightarrow+\infty}
\langle B'(u_n), u_n\rangle=\langle B'(u), u\rangle$;
\item[\rm (iii)] $\lim\limits_{n\rightarrow+\infty}
\langle B'(u_n), \varphi\rangle=\langle B'(u), \varphi\rangle$.
\end{itemize}
\begin{itemize}
\item[\rm (2)] Moreover, for any $u\in H^1(\R^2)$,
\item[\rm (iv)] $$B(u)=\frac{1}{2}\int_{\R^2}u^2\int_{|x|}^\infty\frac{h(s)}{s}u^2(s){\rm d}s{\rm d}x;$$
\item[\rm (v)]  $\lan B'(u), u\ran=6B(u)$.
\end{itemize}
\end{itemize}
\end{lemma}

The following lemma is used to estimate function $h$.
\begin{lemma}\label{Lem:hs}
For any $u\in H_{rad}^1(\R^2)$ and $x\in\R^2$, the following hold:
\begin{itemize}
\item[\rm (i)]  there exists $c_r>0$ such that the $h(|x|)\leq c_r |x|^{\frac{2(r-2)}{r}}\|u\|_r^2$ for $r\in(2,+\infty)$;
\item[\rm (ii)]  for $r\in(2,4)$ and $r'\in(4,+\infty)$, there exist $C_r,C_{r'}>0$ such that
$$
\int_{|x|}^\infty\frac{h(s)}{s}u^2(s){\rm d}s\leq
\|u\|^2_4(C_r\|u\|_r^2+C_r'\|u\|_{r'}^2).
$$ 
\item[\rm (iii)] there holds for $r\in(6,8)$,
$$
\int_{0}^\infty\frac{h(s)}{s}u^2(s){\rm d}s\leq
\left[\frac{2}{\pi}B(u)+\frac{1}{2\pi}\|u\|^2_2+C(\|u\|_{4}^4+\|u\|_r^{\frac{r}{2}})\right].
$$
\end{itemize}
\end{lemma}
\Proof
(i) By H\"{o}lder's inequality, for any $r>2$, there exists $c_r>0$ such that
$$
h(|x|)=\int_{B_{|x|}}\frac{1}{2\pi}u^2(y){\rm d}y
\leq c_r|x|^{\frac{2(r-2)}{r}}\|u\|_r^2.
$$
(ii) It follows from conclusions (i) that if $|x|\leq1$, then
$$
\aligned
\int_{|x|}^\infty\frac{h(s)}{s}u^2(s){\rm d}s
=&\int_{|x|}^{1}\frac{h(s)}{s}u^2(s){\rm d}s+\int_{1}^\infty\frac{h(s)}{s}u^2(s){\rm d}s\\
\leq& C\|u\|^2_4\int_{|x|}^{1}u^2(s){\rm d}s+C\|u\|^2_{r}\int_{1}^\infty s^{\frac{r-4}{r}}u^2(s){\rm d}s\\
\leq& C\|u\|^2_4\left(\int_{|x|}^{1}s^{\frac{-2}{r'-2}}{\rm d}s\right)^{\frac{r'-2}{r'}}
\left(\int_{|x|}^{1}u^{r'}(s)s{\rm d}s\right)^{\frac{2}{r'}}\\
&+
C\|u\|^2_{r}\left(\int_{1}^\infty s^{\frac{r-8}{r}}{\rm d}s\right)^{\frac{1}{2}}\left(\int_{1}^\infty u^4(s)s{\rm d}s\right)^{\frac{1}{2}}\\
\leq& \|u\|^2_4 (C_{r'}\|u\|^2_{r'}+C_r\|u\|^2_{r}).
\endaligned
$$
Based on above, the conclusion holds obviously for the case where $|x|>1$. The proof of conclusion (ii) is complete.\\
(iii) Define
\begin{equation}\label{eqn:K}
\aligned
K_1+K_2&:=\int_{0}^1\frac{h(s)}{s}u^2(s){\rm d}s+\int_{1}^\infty\frac{h(s)}{s}u^2(s){\rm d}s=\int_{0}^\infty\frac{h(s)}{s}u^2(s){\rm d}s.
\endaligned
\end{equation}
One hand, let $\alpha=\frac{r}{8}$, then
\begin{equation}\label{eqn:K1}
\aligned
K_1&=\int_{0}^1\frac{h(s)}{s}u^2(s){\rm d}s+\int_{0}^1\frac{h^\alpha(s)}{s}u^2(s){\rm d}s\\
&\leq\int_{\{s\in(0,1]:\,h(s)\geq1\}}\frac{h(s)}{s}u^2(s){\rm d}s
+\int_{\{s\in(0,1]:\,h(s)\leq1\}}\frac{h(s)}{s}u^2(s){\rm d}s\\
&\leq\int_{0}^1\frac{h(s)}{s}u^2(s){\rm d}s
+\int_{\{s\in(0,1]:\,h(s)\leq1\}}\frac{h^\alpha(s)}{s}u^2(s){\rm d}s\\
&\leq\frac{1}{2\pi}\int_{\R^2}\frac{h^2(|x|)}{|x|^2}u^2(x){\rm d}x
+\int_{0}^1\frac{h^\alpha(s)}{s}u^2(s){\rm d}s\\
&\leq\frac{1}{\pi}B(u)
+\left(\int_{0}^1\frac{h^{2\alpha}(s)}{s^3}{\rm d}s\right)^\frac{1}{2}\|u\|_4^2\\
& \leq\frac{1}{\pi}B(u)+\left(\int_{0}^1s^{\frac{4\alpha(r-2)}{r}-3}\|u\|_r^{4\alpha}{\rm d}s\right)^\frac{1}{2}\|u\|_4^2\\
&\leq\frac{1}{\pi}B(u)+C(\|u\|_4^4+\|u\|_r^{\frac{r}{2}}).
\endaligned
\end{equation}
On the other hand,
\begin{equation}\label{eqn:K2}
\aligned
K_2&=\int_{1}^\infty\frac{h(s)}{s}u^2(s){\rm d}s+\int_{0}^1\frac{h^\alpha(s)}{s}u^2(s){\rm d}s\\
&\leq\int_{\{s\in(1,\infty):\,h(s)\geq1\}}\frac{h(s)}{s}u^2(s){\rm d}s
+\int_{\{s\in(1,\infty):\,h(s)\leq1\}}\frac{h(s)}{s}u^2(s){\rm d}s\\
&\leq\int_{0}^\infty\frac{h^2(s)}{s}u^2(s){\rm d}s
+\int_{1}^\infty u^2(s)s{\rm d}s\\
&=\frac{1}{\pi}B(u)+\frac{1}{2\pi}\|u\|_2^2.
\endaligned
\end{equation}
Substituting (\ref{eqn:K1}) and (\ref{eqn:K2}) into (\ref{eqn:K}), conclusion (iii) holds.
\qed

To prove the existence and multiplicity of nodal solutions,
we recall the following Pohozaev
 identity and an inequality of Sobolev type, which were proved in \cite{Byeon12}.

\begin{lemma}[\cite{Byeon12}]\label{Lem:pohozaev}
Let $b,c$ be positive real constants and $u\in H_{rad}^1(\R^2)$ be weak solution of the equation:
$$
-\triangle u+b u+c\left(\frac{h^2(|x|)}{|x|^2}+\int_{|x|}^{+\infty}\frac{h(s)}{s}u^2(s){\rm d}s\right)u
=\lambda|u|^{p-1}u,\quad x\in\R^2,
$$
 then
$$
b\int_{\R^2}u^2{\rm d}x+2c\int_{\R^2}\frac{h^2(|x|)}{|x|^2}u^2{\rm d}x+\frac{2\lambda}{p}\int_{\R^2}|u|^p{\rm d}x=0.
$$
\end{lemma}
\begin{lemma}[\cite{Byeon12}]\label{Lem:inequality}
For $u\in H_{rad}^1(\R^2)$, there holds
$$
\int_{\R^2}u^4{\rm d}x\leq 4\left(\int_{\R^2}|\nabla u|^2{\rm d}x\right)^\frac{1}{2}\left(\int_{\R^2}\frac{u^2}{|x|^2}
h^2(|x|){\rm d}x\right)^{\frac{1}{2}}.
$$
Furthermore, the equality is attained by a continuum of function
$$
\left\{u_l=\frac{\sqrt{8}l}{1+|lx|^2}\in H_{rad}^1(\R^2)\big|\,l\in(0,\infty)\right\},
$$
and
$$
\frac{1}{4}\int_{\R^2}|u_l|^4{\rm d}x=\int_{\R^2}|\nabla u_l|^2{\rm d}x=\int_{\R^2}\frac{u_l^2}{|x|^2}h_l^2(|x|)
{\rm d}x=\frac{16\pi l^2}{3},
$$
where $h_l(|x|):=\int_{0}^{|x|}\frac{s}{2}u_l^2(s){\rm d}s$.
\end{lemma}


\s{Proof of Theorem \ref{Thm:existence}}
\renewcommand{\theequation}{3.\arabic{equation}}


\subsection{The perturbed problem}
In this section, we investigate
the existence of sign-changing solutions to (\ref{eqn:cs3}) for $p\in(4,6)$
via a perturbation argument.
 In contrast with $p>6$, the {\it Ambrosetti-Rabinowtiz} condition fails for $p\in(4,6)$.
 It results in tough difficulty getting the boundedness of Palais-Smale sequences. In \cite{Byeon12}, the authors constructed a constraint manifold coming from both deformations
of range and domain of function $u$, which was firstly used in \cite{Ruiz06}.
To overcome this difficulty, we introduce a perturbation approach.

For convenience's sake, for any $u\in H_{rad}^1(\R^2)$ and $x\in\R^2\setminus\{0\}$, set
$$
\mathcal{B}(u)(x)=\frac{h^2(|x|)}{|x|^2}+\int_{|x|}^{+\infty}\frac{h(s)}{s}u^2(s){\rm d}s,
$$
then
$$
\lan B'(u),\vp\ran=\int_{\R^2}\mathcal{B}(u)u\vp\,\ud x,\,\,\vp\in H^1(\R^2).
$$
Moreover, by Lemma \ref{Lem:hs}, $\mathcal{B}(u)$ is well defined and $\mathcal{B}(u)(x)\in(0,+\iy),\,\,x\in\R^2\setminus\{0\}$. Taking $$\alpha\in\left(0,\min\left\{\frac{1}{2},\frac{p-4}{2}\right\}\right), \gamma,\beta\in(0,1]\,\, \mbox{and}\,\,\,q>6,$$
we consider the modified problem
\begin{equation}\label{eqn:peturbed}
\aligned
-\triangle u+[\omega+\mathcal{B}(u)]u+\gamma\left(\int_{\R^2}|u|^4{\rm d}x\right)^\alpha u^3=\lambda|u|^{p-2}u+\beta|u|^{q-2}u,\quad x\in\R^2
\endaligned
\end{equation}
and its associated functional is as follows
$$
I_{\gamma,\beta}(u)=I_\lambda(u)+\frac{\gamma}{4(1+\alpha)}\left(\int_{\R^2}|u|^4{\rm d}x\right)^{1+\alpha}-\frac{\beta}{q}\int_{\R^2}|u|^q{\rm d}x.
$$
It is easy to show that $I_{\gamma,\beta}\in C^1(H_{rad}^1(\R^2),\R)$ and
$$
\lan I_{\gamma,\beta}'(u),v\ran= \lan I_\lambda'(u),v\ran +\gamma\|u\|_4^\alpha\int_{\R^2}u^3v{\rm d}x-\beta\int_{\R^2}|u|^{q-2}uv{\rm d}x,\,\,u,v\in H_{rad}^1(\R^2).
$$
For each $u\in H_r^1(\R^2)$, one can show that the equation
\begin{equation} \label{eqn:in1}
\aligned
  -\triangle v+[\omega+\mathcal{B}(u)]v+\gamma\left(\int_{\R^2}|u|^4{\rm d}x\right)^\alpha u^2v=\lambda|u|^{p-2}u+\beta|u|^{q-2}u,\quad x\in\R^2
\endaligned
\end{equation}
has a unique weak solution $v:=T_{\gamma,\beta}(u)\in H_{rad}^1(\R^2)$.
Clearly, the fact that $u$ is a solution of equation (\ref{eqn:in1}) is equivalent to
 that $u$ is a fixed point
of $T_{\gamma,\beta}$.

\begin{lemma}\label{Lem:T1}
The operator $T_{\gamma,\beta}: u\in H_{rad}^1(\R^2)\mapsto v\in H_{rad}^1(\R^2)$ is continuous.
\end{lemma}
\Proof
Assume that $\{u_n\}\subset H_{rad}^1(\R^2)$ with $u_n\rightarrow u$ strongly in $H_{rad}^1(\R^2)$. For simplicity, for any $n$ and $x\in\R^2\setminus\{0\}$, set
$$
\mathcal{B}(u_n)(x)=\frac{h_n^2(|x|)}{|x|^2}+\int_{|x|}^{+\infty}\frac{h_n(s)}{s}u_n^2(s){\rm d}s,
$$
where $h_n(s)=\int_{0}^s\frac{\tau}{2}u_n^2(\tau){\rm d\tau}$.
Let $v=T_{\gamma,\beta}(u)$ and
$v_n=T_{\gamma,\beta}(u_n)$, then we have
\begin{equation}\label{eqn:f-0}
\aligned
&\int_{\R^2}(\nabla v_n\nabla w+[\omega+\mathcal{B}(u_n)]v_nw){\rm d}x
+\gamma\|u_n\|_4^\alpha\int_{\R^2} u_n^2v_nw{\rm d}x
\\&=\lambda\int_{\R^2}|u_n|^{p-2}u_nw{\rm d}x+\beta\int_{\R^2}|u_n|^{q-2}u_nw{\rm d}x, \quad\forall\, w\in H_{rad}^1(\R^2).
\endaligned
\end{equation}
Moreover,
\begin{equation}\label{eqn:f-00}
\aligned
&\int_{\R^2}(\nabla v\nabla w+[\omega+\mathcal{B}(u)]vw){\rm d}x+\gamma\|u\|_4^\alpha\int_{\R^2} u^2vw{\rm d}x
\\&=\lambda\int_{\R^2}|u|^{p-2}uw{\rm d}x+\beta\int_{\R^2}|u|^{q-2}uw{\rm d}x, \quad\forall\, w\in H_{rad}^1(\R^2).
\endaligned
\end{equation}
We show that $\|v_n-v\|\rightarrow0$ as $n\rightarrow\infty$. Indeed, testing with $w=v_n$ in (\ref{eqn:f-0}) gives
$$
\|v_n\|^2\leq\lambda\int_{\R^2}|u_n|^{p-2}u_nv_n{\rm d}x+\beta\int_{\R^2}|u_n|^{q-2}u_nv_n{\rm d}x
$$
which, together with H\"{o}lder's inequality, implies that sequence $\{v_n\}$ is bounded
in $H_{rad}^1(\R^2)$.
Assume that $v_n\rightharpoonup v^*$ weakly in $H_{rad}^1(\R^3)$
and strongly in $L^s(\R^3)$ for $s\in(2,+\infty)$ after extracting a subsequence,
 then by (\ref{eqn:f-0}) and Lemmas \ref{Lem:bu} and \ref{Lem:hs} we have
\begin{equation}\label{eqn:f--}
\aligned
&\int_{\R^3}(\nabla v^*\nabla w+[\omega+\mathcal{B}(u)] v^*w){\rm d}x+
+\gamma\|u\|_4^\alpha\int_{\R^2} u^2v^*w{\rm d}x
\\&=\lambda\int_{\R^2}|u|^{p-2}uw{\rm d}x+\beta\int_{\R^2}|u|^{q-2}uw{\rm d}x \quad\forall w\in H_{rad}^1(\R^2),
\endaligned
\end{equation}
which implies that $v^*$ is a solution of equation (\ref{eqn:in1}). By the uniqueness,
we immediately get $v=v^*$.
Taking $w=v_n-v$ in (\ref{eqn:f-0}) and (\ref{eqn:f-00})
and then subtracting, we have
\begin{align}\label{eqn:f1}
&\|v_n-v\|^2+\int_{\R^2}\mathcal{B}(u_n)(v_n-v)^2{\rm d}x
+\gamma\|u_n\|_4^{\alpha}\int_{\R^2} u_n^2(v_n-v)^2{\rm d}x\nonumber\\
&=\int_{\R^2}[\mathcal{B}(u_n)
-\mathcal{B}(u)] v(v_n-v){\rm d}x+\gamma(\|u_n\|_4^{\alpha}-\|u\|_4^{\alpha})\int_{\R^2}u_n^2v(v_n-v){\rm d}x\\
&\,\,\,\,\,\,\,+\lambda\int_{\R^2}(|u_n|^{p-2}u_n-|u|^{p-2}u)(v_n-v){\rm d}x+\beta\int_{\R^2}(|u_n|^{q-2}u_n-|u|^{q-2}u)(v_n-v){\rm d}x.\nonumber
\end{align}
So it follows from Lemmas \ref{Lem:bu} and \ref{Lem:hs} and Sobolev's embedding inequality that
$v_n\rightarrow v$ in $H_{rad}^1(\R^2)$ as $n\rightarrow \infty$.
Therefore, $T_{\gamma,\beta}$ is continuous.
\qed

\begin{lemma}\label{Lem:T2}
\begin{itemize}
\item[\rm (1) ] $\lan I'_{\gamma,\beta}(u),u-T_{\gamma,\beta}(u)\ran\geq \|u-T_{\gamma,\beta}(u)\|^2$ for all $u\in H_{rad}^1(\R^2)$;
\item[\rm (2) ] $\|I'_{\gamma,\beta}(u)\|\leq C\|u-T_{\gamma,\beta}(u)\|(1+\|u\|^4+\|u\|^2+\|u\|^{2+\alpha})$
for some
$C>0$ and all $u\in H_{rad}^1(\R^2)$.
\end{itemize}
\end{lemma}
\Proof
(1) Since $T_{\lambda,\beta}(u)$ is the solution of equation (\ref{eqn:in1}), we have
$$
\aligned
\lan I'_{\gamma,\beta}(u),u-T_{\gamma,\beta}(u)\ran=&\|u-T_{\gamma,\beta}(u)\|^2
+\int_{\R^2}\mathcal{B}(u)(u-T_{\gamma,\beta}(u))^2{\rm d}x\\
&+\gamma\|u\|_4^{\alpha}\int_{\R^2}u^2|u-T_{\gamma,\beta}(u)|^2{\rm d}x,
\endaligned
$$
which means $\lan I'_{\gamma,\beta}(u),u-T_{\gamma,\beta}(u)\ran\geq \|u-T_{\gamma,\beta}(u)\|^2$
 for all $u\in H_{rad}^1(\R^2)$.\\
 (2) Notice that for all $\varphi\in C_0^\infty(\R^2)\cap H_r^1(\R^2)$,
{\allowdisplaybreaks\begin{align*}
\lan I_{\gamma,\beta}'(u),\varphi\ran=&\int_{\R^2}[\nabla(u-T_{\gamma,\beta}(u))\nabla\varphi
+[\omega+\mathcal{B}(u)](u-T_{\gamma,\beta}(u))\varphi]{\rm d}x\\
&
+\gamma\|u\|_4^{\alpha}\int_{\R^2}u^2(u-T_{\gamma,\beta}(u))\varphi {\rm d}x
\end{align*}}%
which implies by Lemma \ref{Lem:hs} that $\|I_{\gamma,\beta}'(u)\|\leq C\|u-T_{\gamma,\beta}(u)\|(1+\|u\|^4
+\|u\|^{2}+\|u\|^{2+\alpha})$ for all $u\in H_{rad}^1(\R^2)$.
\qed
\begin{lemma}\label{Lem:T3}
For fixed $(\gamma,\beta)\in(0,1]\times(0,1]$, $c<d$ and $a>0$, there exists $\varepsilon_{\gamma,\beta}>0$ such that
$\|u-T_{\gamma,\beta}(u)\|\geq \varepsilon_{\gamma,\beta}$
if $u\in H_{rad}^1(\R^2)$, $I_{\gamma,\beta}(u)\in [c,d]$
and $\|I_{\gamma,\beta}'(u)\|\geq a$.
\end{lemma}
\Proof
Fix $\mu\in(6,q)$, then for all $u\in H_{rad}^1(\R^2)$, we have
\begin{equation}\label{eqn:T0}
\aligned
&I_{\gamma,\beta}(u)-\frac{1}{\mu}\langle u,u-T_{\gamma,\beta}(u)\rangle\\
&=\frac{\mu-2}{2\mu}\|u\|^2
+\frac{1}{\mu}\int_{\R^2}\mathcal{B}(u)u(u-T_{\gamma,\beta}(u)){\rm d}x
+\frac{\mu-2}{\mu}B(u)\\
&\,\,\,\,\,+\frac{\gamma}{\mu}\|u\|_4^{\alpha}\int_{\R^2}u^3(u-T_{\gamma,\beta}(u)){\rm d}x-\frac{1}{\mu}\int_{\R^2}u^2\int_{|x|}^{+\infty}\frac{h(s)}{s}u^2(s){\rm d}s{\rm d}x
\\
&\,\,\,\,\,+\frac{(p-\mu)\lambda}{\mu p}\int_{\R^2}|u|^{p}{\rm d}x+\frac{(q-\mu)\beta}{\mu q}\int_{\R^2}|u|^{q}{\rm d}x+\gamma\frac{\mu-4}{4\mu}\|u\|_4^{4(1+\alpha)}.
\endaligned
\end{equation}
Then, by the conclusion (iv) of Lemma \ref{Lem:bu}, we have
\begin{equation}\label{eqn:T2}
\aligned
&\|u\|^2+2B(u)+\beta\|u\|_q^q+\gamma\|u\|_4^{4(1+\alpha)}-\lambda\|u\|_{p}^{p}\\
&\le c\left\{|I_{\gamma,\beta}(u)|+\|u\|\|u- T_{\gamma,\beta}(u)\|
+\int_{\R^2}\mathcal{B}(u)|u(u-T_{\gamma,\beta}(u))|{\rm d}x\right.\\
&\,\,\,\,\,\,\,\,\,\,\,\,\,\,\left.+\frac{\gamma}{\mu}\|u\|_4^{\alpha}\int_{\R^2}|u^3(u-T_{\gamma,\beta}(u))|{\rm d}x\right\},
\endaligned
\end{equation}
where $c$ depends on $\beta,\gamma,\mu,\alpha,\lambda$.
One hand, take $r\in(\frac{4(1+\alpha)}{1+2\alpha},4)$ and $r'=q$ in Lemma \ref{Lem:hs}, then
it follows from Lemma \ref{Lem:hs}, Jensen's inequality,
H\"{o}lder's inequality and Sobolev's inequality, we have for any $\xi>0$,
there exists $C_\xi>0$ such that
\begin{equation}\label{eqn:T3-0}
\aligned
&\int_{\R^2}|u(u-T_{\gamma,\beta}(u))|\int_{|x|}^{+\infty}\frac{h(s)}{s}u^2(s){\rm d}s{\rm d}x\\
&\leq\left(\int_{\R^2}u^2\int_{|x|}^{+\infty}\frac{h(s)}{s}u^2(s){\rm d}s{\rm d}x\right)^\frac{1}{2}
\left(\int_{\R^2}(u-T_{\gamma,\beta}(u))^2\int_{|x|}^{+\infty}\frac{h(s)}{s}u^2(s){\rm d}s{\rm d}x\right)^\frac{1}{2}\\
&\leq\xi B(u)+C_\xi \|u\|_4^2\|u\|_r^2\|u-T_{\gamma,\beta}(u)\|^2+C_\xi \|u\|_4^2\|u\|_{q}^2\|u-T_{\gamma,\beta}(u)\|^2\\
&\leq\xi B(u)+C_\xi \left(\|u\|_4^{4(1+\alpha)}\|u-T_{\gamma,\beta}(u)\|^{1+\alpha}+
\|u\|_r^r
\|u-T_{\gamma,\beta}(u)\|^{\frac{r}{4}}\right.\\
&\,\,\,\,\,\,\,\,\,\,\,\,\,\,\,\,\,\,\,\,\,\,\,\,\,\,\,\,\,\,\,\,\,\,\,\,\,\,\,\,\,\,\left.+\|u-T_{\gamma,\beta}(u)\|^{\frac{2r(1+\alpha)}{(1+2\alpha)r-4(1+\alpha)}}+(\|u\|_4^4+\|u\|_q^4)\|u-T_{\gamma,\beta}(u)\|^2\right).
\endaligned
\end{equation}
On the other hand, it follows from
H\"{o}lder's inequality and Sobolev's inequality, we have for any $\xi>0$,
there exists $C_\xi>0$ such that
\begin{equation}\label{eqn:T3}
\aligned
&\int_{\R^2}\frac{h^2(|x|)}{|x|^2}u(u-T_{\gamma,\beta}(u)){\rm d}x\\
&\leq\left(\int_{\R^2}\frac{h^2(|x|)}{|x|^2}u^2{\rm d}x\right)^\frac{1}{2}
\left(\int_{\R^2}\frac{h^2(|x|)}{|x|^2}(u-T_{\gamma,\beta}(u))^2{\rm d}x\right)^\frac{1}{2}\\
&\leq\xi B(u)+C_\xi \|u\|_4^4\|u-T_{\gamma,\beta}(u)\|^2.
\endaligned
\end{equation}
Combining (\ref{eqn:T2})-(\ref{eqn:T3}) we get
\begin{equation}\label{eqn:T4}
\aligned
&\|u\|^2+ B(u)+\beta\|u\|_q^q+\gamma\|u\|_4^{4(1+\alpha)}-\lambda\|u\|_{p}^{p}\\
&\le C\left(|I_{\gamma,\beta}(u)|+\|u\|\|u- T_{\gamma,\beta}(u)\|+\|u\|_r^r
\|u-T_{\gamma,\beta}(u)\|^{\frac{r}{4}}\right.\\
&\,\,\,\,\,\,\,\,\,\,\,\,\ \left.+\|u\|_4^{4(1+\alpha)}(\|u- T_{\gamma,\beta}(u)\|^4+\|u-T_{\gamma,\beta}(u)\|^{1+\alpha})\right.\\
&\,\,\,\,\,\,\,\,\,\,\,\,\,\, \left.+\|u-T_{\gamma,\beta}(u)\|^{\frac{2r(1+\alpha)}{(1+2\alpha)r-4(1+\alpha)}}+(\|u\|_4^4+\|u\|_q^4)\|u-T_{\gamma,\beta}(u)\|^2\right).
\endaligned
\end{equation}
Assume on the contrary that there exists $\{u_n\}\subset H_{rad}^1(\R^2)$
with $I_{\gamma,\beta}(u_n)\in[c,d]$ and
$\|I'_{\gamma,\beta}(u_n)\|\geq \alpha$ such that
$$
\|u_n-T_{\gamma,\beta}(u_n)\|\rightarrow0,\quad\text{as}\,\,n\rightarrow\infty.
$$
Then owing to $r\in(2,q)$, we have
\begin{equation}\label{eqn:T4-1}
\int_{\R^2}\left(\frac{\omega}{2}|u_n|^2+\frac{\beta}{2}|u_n|^q-C\|u_n-T_{\gamma,\beta}(u_n)\|^{\frac{r}{2}}|u_n|^r\right){\rm d}x\geq0
\end{equation}
for $n$ large enough, where $C$ has appeared in (\ref{eqn:T4}). Now we claim that sequence $\{u_n\}$ is bounded in $H_{rad}^1(\R^2)$.
Otherwise, assume $\|u_n\|\rightarrow\infty$.
From  (\ref{eqn:T4}) and (\ref{eqn:T4-1}), we deduce that
\begin{equation}\label{eqn:T5}
\|u_n\|^2+ B(u_n)+\beta\|u_n\|_q^q+\gamma\|u_n\|_4^{4(1+\alpha)}
-\lambda\|u_n\|_{p}^{p}\leq C
\end{equation}
for large $n$. Here $$
B(u_n):=\frac{1}{2}\int_{\R^2}\frac{u_n^2}{|x|^2}h_n^2(|x|){\rm d}x.
$$
Note that, for any $A_1>0$, we can choose $A_2>0$ such that for any $t\in\R$,
$$
t^{1+\alpha}>A_1t-A_2.
$$
Applying this with $t=\|u_n\|_4^4$, then by (\ref{eqn:T5}) we have
\begin{equation}\label{eqn:T7}
\|u_n\|^2+B(u_n)
+\int_{\R^2}(\beta|u_n|^q+\gamma A_1|u_n|^{4}-\lambda|u_n|^p){\rm d}x-A_2\leq C.
\end{equation}
Since $4<p<q$, we take $A_1$ large enough such that the function
$\gamma A_1|t|^{4}+\beta|t|^q-C |t|^p>0$ for any $t\in\R$.
So it is easy to see from (\ref{eqn:T7}) that $\{u_n\}$ is bounded
in $H_{rad}^1(\R^2)$ for any fixed $(\gamma,\beta)\in(0,1]\times(0,1]$.
The claim combined with Lemma \ref{Lem:T2} implies $\|I'_{\gamma,\beta}(u_n)\|\rightarrow0$
as $n\rightarrow\infty$, which is a contradiction.
\qed

\subsection{Invariant Subsets of Descending Flows}

Define the positive and negative cones by
$$
P^+:=\{u\in H_{rad}^1(\R^2):\, u\geq0\}\quad\text{and}\quad P^-:=\{u\in H_{rad}^1(\R^2):\,u\leq0\},
$$
respectively. Set for $\epsilon>0$,
$$
P_\epsilon^+:=\{u\in H_{rad}^1(\R^2):\, dist(u,P^+)<\epsilon\}\quad\text{and}\quad P_\epsilon^-:=\{u\in H_{rad}^1(\R^2):\,dist(u,P^-)<\epsilon\},
$$
where $dist(u,P^\pm)=\inf\limits_{v\in P^\pm}\|u-v\|$. Clearly, $P_\epsilon^-=-P_\epsilon^+$. Let $W=P_\epsilon^+\cup P_\epsilon^-$.
It is not hard to check that $W$ is an open and symmetric subset of $H_{rad}^1(\R^2)$
and $H_{rad}^1(\R^2)\setminus {W}$ contains only sign-changing functions.

The following will prove that, for $\epsilon$ small enough,
all sign-changing solutions to equation (\ref{eqn:peturbed}) are contained in
$H_{rad}^1(\R^2)\setminus {W}$.

\begin{lemma}\label{Lem:C1}
 There exists $\epsilon_0>0$ such that for $\epsilon\in(0,\epsilon_0)$,
\item[\rm (1) ] $T_{\gamma,\beta}(\partial P_\epsilon^-)\subset P_\epsilon^-$ and every nontrivial solution $u\in P_\epsilon^-$ is negative,
\item[\rm (2) ] $T_{\gamma,\beta}(\partial P_\epsilon^+)\subset P_\epsilon^+$ and every nontrivial solution $u\in P_\epsilon^+$ is positive.
\end{lemma}
\Proof
Since the two conclusions are similar, we only prove $T_{\gamma,\beta}(\partial P_\epsilon^-)\subset P_\epsilon^-$.
Let $u\in H_{rad}^1(\R^2)$ and $v=T_{\gamma,\beta}(u)$. Since $\text{dist}(v,P^-)\leq \|v^+\|$, by Sobolev's inequality
we have
$$
\aligned
&\text{dist}(v,P^-)\|v^+\|=\|v^+\|^2=\langle v,v^+\rangle\\
&\leq\lambda\int_{\R^2}|u|^{p-2}uv^+{\rm d}x+\beta\int_{\R^2}|u|^{q-2}uv^+{\rm d}x
-\gamma\|u\|_4^{\alpha}\int_{\R^2}u^2vv^+{\rm d}x\\
&\le\lambda\int_{\R^2}|u^+|^{p-2}u^+v^+{\rm d}x+\beta\int_{\R^2}|u^+|^{q-2}u^+v^+{\rm d}x\\
&\leq C[\text{dist}(u,P^-)^{p-1}+\text{dist}(u,P^-)^{q-1}]\|v^+\|,
\endaligned
$$
which implies that
$$
\text{dist}(v,P^-)\le C[\text{dist}(u,P^-)^{p-1}+\text{dist}(u,P^-)^{q-1}].
$$
Then there exists $\epsilon_0>0$ such that for $\epsilon\in(0,\epsilon_0)$,
$$
\text{dist}(T_{\gamma,\beta}(u),P^-)=\text{dist}(v,P^-)<\epsilon.
$$
Therefore, we have $T_{\gamma,\beta}(u)\in P_\epsilon^-$ for any $u\in P_\epsilon^-$.
Assume that there exists $u\in P_\epsilon^-$ such that $T_{\gamma,\beta}(u)=u$, then
$u\in P^-$.\ If $u\not\equiv 0$, by the maximum principle, $u<0$ in
$\R^2$.
\qed

Since the operator $T_{\gamma,\beta}$ may not be locally Lipschitz continuous,
we need to construct a locally Lipschitz continuous vector field which inherits its properties.
Similar to the proof of Lemma 2.1 in \cite{Bartsch05}, we have
\begin{lemma}\label{Lem:local}
There exists a locally Lipschitz continuous operator
$B_{\gamma,\beta}:\,H_{rad}^1(\R^2)\setminus K_{\gamma,\beta}\rightarrow H_{rad}^1(\R^2)$ such that
\begin{itemize}
\item[(i)] $\langle T'_{\gamma,\beta}(u),u-B_{\gamma,\beta}(u)\rangle\geq\frac{1}{2}\|u-T_{\gamma,\beta}(u)\|^2$;
\item[ (ii)]$\frac{1}{2}\|u-B_{\gamma,\beta}(u)\|^2\leq\|u-T_{\gamma,\beta}(u)\|^2\leq2
\|u-B_{\gamma,\beta}(u)\|^2$;
\item[ (iii)]$T_{\gamma,\beta}(\partial P_\epsilon^\pm)\subset P_\epsilon^\pm$, $\forall \epsilon\in(0,\epsilon_0)$;
\item[ (iv)] $B_{\gamma,\beta}$ is odd,
\end{itemize}
where $K_{\gamma,\beta}:=\{u\in H_{rad}^1(\R^2):\,I_{\gamma,\beta}'(u)=0\}$
and $\epsilon_0$ has been given in Lemma \ref{Lem:C1}.
\end{lemma}
In what follows, we prove functional $I_{\gamma,\beta}$ satisfies the (PS)-condition.

\begin{lemma}\label{Lem:PS}
Let $(\gamma,\beta)\in(0,1)\times(0,1)$ and $c\in\R$. Assume $\{u_n\}\subset H_{rad}^1(\R^2)$ satisfy $I_{\gamma,\beta}(u_n)\rightarrow c$ and
$I'_{\gamma,\beta}(u_n)\rightarrow 0$  as $n\rightarrow\infty$,
then up to a subsequence, $u_n\rightarrow u$ in $H_r^1(\R^2)$ for some $u\in H_{rad}^1(\R^2)$.

\end{lemma}
\Proof
For $\mu\in(6,q)$, we have

$$
\aligned
&I_{\gamma,\beta}(u_n)-\frac{1}{\mu}\lan I'_{\gamma,\beta}(u_n),u_n\ran=\frac{\mu-2}{2\mu}\|u_n\|^2+\frac{\mu-6}{\mu}B(u_n)\\
&+\gamma\frac{\mu-4(1+\alpha)}{4(1+\alpha)\mu}\|u_n\|_4^{4(1+\alpha)}
-\lambda\frac{\mu-p}{p\mu}\int_{\R^2}|u_n|^p{\rm d}x+\beta\frac{q-\mu}{q\mu}\int_{\R^2}|u|^q{\rm d}x,
\endaligned
$$
which implies that
$$
\aligned
& |I_{\gamma,\beta}(u_n)|+\frac{1}{\mu}\|I'_{\gamma,\beta}(u_n)\|\|u_n\|\geq\frac{\mu-2}{2\mu}\|u_n\|^2+\frac{\mu-6}{\mu}B(u_n)\\
&+\gamma\frac{\mu-4(1+\alpha)\mu}{4(1+\alpha)}\|u_n\|_4^{4(1+\alpha)}
-\lambda\frac{\mu-p}{p\mu}\int_{\R^2}|u_n|^p{\rm d}x+\beta\frac{q-\mu}{q\mu}\int_{\R^2}|u|^q{\rm d}x.
\endaligned
$$
As in the proof of Lemma \ref{Lem:T3}, one sees that sequence $\{u_n\}$ is bounded in $H_{rad}^1(\R^2)$.
Without loss of generality, we assume that there exists $u\in H_{rad}^1(\R^2)$ such that
$$
\aligned
&u_n\rightharpoonup u\,\,\text{weakly\,in}\,\,H_{rad}^1(\R^2),\\
& u_n\rightarrow u\,\,\text{strongly\,in}\,\,L^r(\R^2)\,\text{for}\,r\in(2,+\infty).
\endaligned
$$
Note that
\begin{equation}\label{eqn:PS1}
\aligned
&\langle I'_{\gamma,\beta}(u_n)-I'_{\gamma,\beta}(u),u_n-u\rangle\\
&=\|u_n-u\|^2+\int_{\R^2}\mathcal{B}(u_n)(u_n-u)^2{\rm d}x+\int_{\R^2}[\mathcal{B}(u_n)-\mathcal{B}(u)]u(u_n-u){\rm d}x\\
&\,\,\,\,\,\,\,-\lambda\int_{\R^2}(|u_n|^{p-2}u_n-|u|^{p-2}u)(u_n-u){\rm d}x+\gamma\|u_n\|_4^{\alpha}\int_{\R^2}u_n^2(u_n-u)^2{\rm d}x\\
&\,\,\,\,\,\,\,+\gamma(\|u_n\|_4^{\alpha}-\|u\|_4^{\alpha})\int_{\R^2}(|u_n|^{2}u-|u|^{2}u)(u_n-u){\rm d}x\\
&\,\,\,\,\,\,\,-\beta\int_{\R^2}(|u_n|^{q-2}u_n-|u|^{q-2}u)(u_n-u){\rm d}x.
\endaligned
\end{equation}
Based on the boundedness of $\{u_n\}$ in $H_{rad}^1(\R^2)$ and the conclusions of Lemma \ref{Lem:hs} we can deduce that
$$
\aligned
&\int_{\R^2}\mathcal{B}(u_n)(u_n-u)^2{\rm d}x\rightarrow0,\quad\int_{\R^2}[\mathcal{B}(u_n)-\mathcal{B}(u)]u(u_n-u){\rm d}x\rightarrow0,\\
&\int_{\R^2}(|u_n|^{p-2}u_n-|u|^{p-2}u)(u_n-u){\rm d}x\rightarrow0,\quad\langle I'_{\gamma,\beta}(u_n)-I'_{\gamma,\beta}(u),u_n-u\rangle\rightarrow0,\\
&\|u_n\|_4^{\alpha}\int_{\R^2}u_n^2(u_n-u)^2{\rm d}x\rightarrow0,\quad(\|u_n\|_4^{\alpha}-\|u\|_4^{\alpha})\int_{\R^2}(|u_n|^{2}u-|u|^{2}u)(u_n-u){\rm d}x\rightarrow0,\\
&\int_{\R^2}(|u_n|^{q-2}u_n-|u|^{q-2}u)(u_n-u){\rm d}x\rightarrow0,
\quad \text{as}\,\,n\rightarrow\infty.
\endaligned
$$
Therefore, from (\ref{eqn:PS1}) we deduce that $u_n\rightarrow u$ in $H_{rad}^1(\R^2)$ as $n\rightarrow\infty$.
\qed
\vskip0.1in
Now we recall a critical point theorem introduced in \cite{Liuj15}. Let $E$ be a Banach space,
$J\in C^1(E,\R)$.
we refer to
$P,Q\subset E$ be open sets, $M=P\cap Q$,
$\Sigma=\partial P\cap\partial Q$ and $W=P\cup Q$.
We denote by $K$ the set of critical point of $J$, that is,
$K=\{u\in E:\,\,J'(x)=0\}$ and $E_0:=E\setminus{K}$.
For $c\in\R$, $K_c=\{x\in E:J(x)=c,J'(x)=0\}$ and
$J^c=\{x\in E:J(x)\leq c\}$.

\begin{definition}{\rm\cite{Liuj15}}\label{Def:1}
$\{P,Q\}$ is called an admissible family of invariant sets
with respect to $J$ at level $c$, provided that the following deformation property
holds: if $K_c\setminus{W }=\emptyset$, then, there exists $\epsilon_0 > 0$ such that for
$\epsilon\in (0,\epsilon_0)$, there
exists $\eta\in C(E,E)$ satisfying
\begin{itemize}
\item[\rm (1) ]
$\eta(\bar{P})\subset\bar{P}$, $\eta(\bar{Q})\subset\bar{Q}$;
\item[\rm (2) ]$\eta|_{J^{c-2\epsilon}}=id$;
\item[\rm (3) ] $\eta(J^{c+\epsilon}\setminus{W})\subset J^{c-\epsilon}$.
\end{itemize}
\end{definition}

\begin{theorem}{\rm\cite{Liuj15}}\label{Thm:xjl}
 Assume that $\{P,Q\}$ is an admissible family of invariant
sets with respect to $J$ at any level $c\geq c_* := \inf_{u\in\Sigma}J(u)$ and there exists a map
$\psi_0 :\triangle\rightarrow E$ satisfying
\begin{itemize}
\item[\rm (1) ]
$\psi_0(\partial_1\triangle)\subset P$ and $\psi_0(\partial_2\triangle)\subset Q$,
\item[\rm (2) ]$\psi_0(\partial_0\triangle)\cap M=\emptyset$,
\item[\rm (3) ] $\sup\limits_{u\in\psi_0(\partial_0\triangle)}J(u)<c_*$,
\end{itemize}
where $\triangle=\{(t_1,t_2)\in\R^2:t_1,t_2>0,\,t_1+t_2\leq1\}$, $\partial_1\triangle=\{0\}\times[0,1]$,
$\partial_2\triangle=[0,1]\times\{0\}$ and  $\partial_0\triangle=\{(t_1,t_2)\in\R^2:t_1,t_2\geq0,\,t_1+t_2=1\}$.
Define
$$
c=\inf\limits_{\psi\in\Gamma}\sup\limits_{u\in\psi(\triangle)\setminus{W}}J(u),
$$
where $\Gamma:=\{\psi\in C(\triangle,E):\,\psi(\partial_1\triangle)\subset P,\,
\psi(\partial_2\triangle)\subset Q,\,\psi|_{\partial_0\triangle}=\psi_0|_{\partial_0\triangle}\}$.
Then $c\geq c_*$ and $K_c\setminus{W}\not=\emptyset$.
\end{theorem}
In order to use Theorem \ref{Thm:xjl} to prove the existence of sign-changing solutions
to equation (\ref{eqn:peturbed}), we take $E=H_{rad}^1(\R^2)$,
$P=P_\epsilon^+,Q=P_\epsilon^-$ and $J=I_{\gamma,\beta}$. We need to show that $\{P_\epsilon^+,P_\epsilon^-\}$
is an admissible family of invariant sets for the
functional $I_{\lambda,\beta}$ at any level $c\in \R$.
Indeed, $K_c\subset W$ if $K_c \setminus{W}=\emptyset$. Since
the functional $I_{\gamma,\beta}$ satisfies the (PS)-condition, $K_c$ is compact. Thus, $2\delta:= dist(K_c; \partial W)>0$.
Similar to Lemma 3.6 in \cite{Liuz16}, we have
\begin{lemma}\label{Lem:XB}(Deformation lemma)
Let $S\subset H_{rad}^1(\R^2)$ and $c\in\R$ such that
$$
\forall u\in I_{\gamma,\beta}^{-1}([c-2\epsilon_0,c+2\epsilon_0])\cap S_{2\delta}:\,\,\,\|I'_{\gamma,\beta}(u)\|\geq \epsilon_0,
$$
where $\epsilon_0$ has been given in Lemma \ref{Lem:C1} and $S_{2\delta}:=\{u\in S:\,\,\text{dist}(u,S)<2\delta\}$. Then for
$\epsilon_1\in(0,\epsilon_0)$ there exists $\eta\in C([0,1]\times H_{rad}^1(\R^2),H_{rad}^1(\R^2))$ such that
\begin{itemize}
\item[(i)] $\eta(t,u)=u$, if $t=0$ or if $u\not\in I^{-1}_{\gamma,\beta}([c-2\epsilon_1,c+2\epsilon_1])$;
\item[ (ii)]$\eta(1,I^{c+\epsilon_1}_{\gamma,\beta}\cap S)\subset I^{c-\epsilon_1}_{\gamma,\beta}$;
\item[ (iii)] $I_{\gamma,\beta}(\eta(\cdot,u))$ is not increasing for all $u\in H_{rad}^1(\R^2)$;
\item[ (iv)]$\eta(t, \overline{P_\epsilon^+})\subset \overline{P_\epsilon^+}$, $\eta(t,\overline{P_\epsilon^-})\subset \overline{P_\epsilon^-}$, $\forall t\in [0,1]$.
\end{itemize}
\end{lemma}

\subsection{Proof of Theorem \ref{Thm:existence}}

\begin{lemma}\label{Lem:qfanshu}
For $r\in[2,+\infty)$ there exits $C>0$ independent of $\epsilon$ such that
$\|u\|_r\leq C\epsilon$ for $u\in M=P_\epsilon^+\cap P_\epsilon^-$.
\end{lemma}
\Proof
For any fixed $u\in M$, we have
$$
\|u^\pm\|_r=\inf\limits_{v\in P^{\mp}}\|u-v\|_r\leq C\inf\limits_{v\in P^\mp}\|u-v\|\leq C\text{dist}(u,P^{\mp}).
$$
So $\|u\|_r\leq C\epsilon$ for $u\in M$.
\qed

\begin{lemma}\label{Lem:xiaozhi}
$I_{\gamma,\beta}(u)\geq\frac{\epsilon^2}{4}$ for all $u\in \partial P_\epsilon^+\cap\partial P_\epsilon^-$ if $\epsilon>0$ is small enough,
that is, $c_*\geq\frac{\epsilon^2}{4}$.
\end{lemma}
\Proof
For any fixed $u\in \partial P_\epsilon^+\cap\partial P_\epsilon^-$, we have $\|u^+\|\geq\text{dist}(u,P^{\mp})=\epsilon$.
By Lemma \ref{Lem:qfanshu}, we have
$$
\aligned
I_{\gamma,\beta}(u)&=\frac{1}{2}\|u\|^2+B(u)
+\frac{\gamma}{4(1+\alpha)}\|u\|_4^{4(1+\alpha)}-\int_{\R^2}\left(\frac{\lambda}{p}|u|^p+\frac{\beta}{q}|u|^q\right){\rm d}x\\
&\geq\frac{1}{2}\epsilon^2-\frac{\lambda}{p}\|u\|_p^p-\frac{1}{q}\|u\|_q^q\\
&\geq\frac{1}{2}\epsilon^2-\frac{C}{p}\epsilon^p-\frac{C}{q}\epsilon^q\geq\frac{\epsilon^2}{4}
\endaligned
$$
for $\epsilon$ small enough.
\qed
\vskip0.1in

\textbf{Proof of Theorem 1.1.} We use Theorem \ref{Thm:xjl} to prove the existence of sign-changing
solutions to equation (\ref{eqn:peturbed}).
Take $S=H_{rad}^1(\R^2)\setminus{W}$ in Lemma \ref{Lem:XB}, then we can easily deduce that $\{P_\epsilon^+,P_\epsilon^-\}$ is an admissible
family of invariant sets for the functional $I_{\gamma,\beta}$ at any level $c\in\R$.\\
In what follows, we divide three steps to complete the proof.\\
\textbf{Step 1}. Choose $v_1,v_2\in C_0^{\infty}(B_1(0))$ such that $\text{supp}(v_1)\cap\text{supp}(v_2)=\emptyset$ and
$v_1<0,v_2>0$, where $B_r(0):=\{x\in\R^2:\,|x|<r\}$. For $(t,s)\in \triangle$, let
$$
\psi_0(t,s):=R[tv_1(R\cdot)+sv_2(R\cdot)],
$$
where $R>1$ will be determined later.
Clearly, for $t,s\in[0,1]$, $\psi_0(0,s)(\cdot)=Rsv_2(R\cdot)\in P_\epsilon^+$
and $\psi_0(t,0)(\cdot)=Rtv_1(R\cdot)\in P_\epsilon^-$.
In view of Lemma \ref{Lem:xiaozhi}, we have for small $\epsilon>0$,
 $c_*=\inf_{u\in\Sigma}I_{\gamma,\beta}(u)\geq\frac{\epsilon^2}{4}$
for any $(\lambda,\beta)\in(0,1]\times(0,1]$.
Set
$$
\rho:=\min\{\|tv_1+(1-t)v_2\|_2:\,\,0\leq t\leq 1\}>0.
$$
then $\|u_t\|_2^2\geq\rho$ for $u\in\psi_0(\partial_0\triangle)$ and
it follows from Lemma \ref{Lem:qfanshu} that
$\psi_0(\partial_0\triangle)\cap P_\epsilon^+\cap P_\epsilon^-=\emptyset$.
Let $u_t=\psi_0(t,1-t)$ for $t\in [0,1]$.  A direct computation shows that
$$
\aligned
&\int_{\R^2}|\nabla u_t|^2{\rm d}x=R^{2}\int_{\R^2}(t^2|\nabla v_1|^2+(1-t)^2|\nabla v_2|^2){\rm d}x,\\
&\int_{\R^2}|u_t|^r{\rm d}x=R^{r-2}\int_{\R^2}(t^r|v_1|^r+(1-t)^r|v_2|^r){\rm d}x\,\,\text{for}\,\,r\in[2,+\infty),\\
&\left(\int_{\R^2}|u_t|^4{\rm d}x\right)^{1+\alpha}=R^{2(1+\alpha)}
\left(\int_{\R^2}(t^4|v_1|^4+(1-t)^4|v_2|^4){\rm d}x\right)^{1+\alpha},\\
& B(u_t)\le R^{2}B(\bar{u}_t),\quad \text{where}\,\,\bar{u}_t=tv_1+(1-t)v_2.
\endaligned
$$
Based on the above facts, we have
$$
\aligned
I_{\gamma,\beta}(u_t)\leq& \frac{1}{2}\|u_t\|^2+B(u_t)
+\frac{\gamma}{4(1+\alpha)}\left(\int_{\R^2}|u_t|^4{\rm d}x\right)^{1+\alpha}
-\int_{\R^2}\left(\frac{\lambda}{p}|u_t|^p+\frac{\beta}{q}|u_t|^q\right){\rm d}x\\
\leq& \frac{R^{2}}{2}\int_{\R^2}(t^2|\nabla v_1|^2+(1-t)^2|\nabla v_2|^2){\rm d}x
+R^{2}B(\bar{u}_t)\\
&+\frac{\gamma R^{2(1+\alpha)}}{4(1+\alpha)}
\left(\int_{\R^2}(t^4|v_1|^4+(1-t)^4|v_2|^4){\rm d}x\right)^{1+\alpha}\\
&+\frac{1}{2}\int_{\R^2}\omega(t^2|v_1|^2+(1-t)^2|v_2|^2){\rm d}x
-\lambda \frac{R^{p-2}}{p}\int_{\R^2}(t^p|v_1|^p+(1-t)^p|v_2|^p){\rm d}x.
\endaligned
$$
Due to $0<\alpha<\frac{p-4}{2}$,
one sees that $I_{\gamma,\beta}(u_t)\rightarrow -\infty$ as $R\rightarrow\infty$ uniformly
for $(\gamma,\beta)\in(0,1)\times(0,1), t\in[0,1]$.
Hence, choosing $R$ large enough, we have
$$
\sup\limits_{u\in\psi_0(\partial_0\triangle)}I_{\gamma,\beta}(u)
<c_*:=\inf\limits_{u\in\Sigma}I_{\gamma,\beta}(u),\quad\forall (\gamma,\beta)\in(0,1)\times(0,1).
$$
Now we have known that $I_{\gamma,\beta}$ satisfies the assumptions of Theorem \ref{Thm:xjl}.
Therefore, the number
$$
c_{\gamma,\beta}=\inf\limits_{\psi\in\Gamma}\sup\limits_{u\in\psi(\triangle)
\setminus{W}}I_{\gamma,\beta}(u)
$$
is a critical value of $I_{\gamma,\beta}$ satisfying $c_{\gamma,\beta}\geq c_*$ and
there exists $u_{\gamma,\beta}\in H_{rad}^1(\R^2)\setminus{(P_\epsilon^+\cup P_\epsilon^-)}$ such that
$I_{\gamma,\beta}(u_{\gamma,\beta})=c_{\gamma,\beta}$ and $I'_{\gamma,\beta}(u_{\gamma,\beta})=0$
for $(\gamma,\beta)\in(0,1]\times(0,1]$.

\textbf{Step 2.}
Passing to the limit as $\gamma\rightarrow0$ and $\beta\rightarrow0$. From the definition of $c_{\gamma,\beta}$,
we see that for any $(\gamma,\beta)\in(0,1]\times(0,1]$,
\begin{equation}\label{eqn:fun0}
c_{\gamma,\beta}\leq C_R:=\sup\limits_{u\in\psi_0(\triangle)}I_{1,0}(u)<\infty,
\end{equation}
where $C_R$ is independent of $(\gamma,\beta)\in(0,1]\times(0,1]$.
Without loss of generality, we set $\gamma=\beta$. Choosing a sequence $\{\gamma_n\}\subset(0,1]$
such that $\gamma_n\rightarrow0^+$, then
we find a sequence of sign-changing critical points $\{u_{\gamma_n}\}$(still denoted by $\{u_n\}$) of $I_{\gamma_n,\beta_n}$ and
$\lim_{n\rightarrow\infty}I_{\gamma_n,\beta_n}(u_n)=\lim_{n\rightarrow\infty}c_{\gamma_n,\beta_n}
=c^*\geq c_*$.
Now we show that $\{u_n\}$ is a bounded sequence in $H_{rad}^1(\R^2)$.
 By the definition of $I_{\gamma,\beta}$, we have
\begin{equation}\label{eqn:fun1}
\aligned
c_{\gamma_n,\beta_n}=&\frac{1}{2}\|u_n\|^2+B(u_n)+\frac{\gamma}{4(1+\alpha)}\|u_n\|_4^{4(1+\alpha)}-\int_{\R^2}\left(\frac{\lambda}{p}|u_n|^p+\frac{\beta}{q}|u_n|^q\right){\rm d}x
\endaligned
\end{equation}
and
 \begin{equation}\label{eqn:fun2}
\aligned
\|u_n\|^2+6B(u_n)
+\gamma\|u_n\|_4^{4(1+\alpha)}-\int_{\R^2}\left(\lambda|u_n|^p+\beta|u_n|^q\right){\rm d}x=0.
\endaligned
\end{equation}

Moreover, it follows from the Pohozaev identity that
 \begin{equation}\label{eqn:fun3}
\aligned
\omega \|u_n\|_2^2
+4B(u_n)+\frac{\gamma}{2}\|u_n\|_4^{4(1+\alpha)}-\int_{\R^2}\left(\frac{2\lambda}{p}|u_n|^p+\frac{2\beta}{q}|u_n|^q\right){\rm d}x=0.
\endaligned
\end{equation}
Multiplying (\ref{eqn:fun1}) by $2(1+\alpha)$, (\ref{eqn:fun2}) by $-1$ and (\ref{eqn:fun3}) by $1$ and adding them up,
we get
\begin{equation}\label{eqn:fun4}
\aligned
2(1+\alpha)c_{\gamma_n,\beta_n}=&\alpha\int_{\R^2}|\nabla u_n|^2{\rm d}x+(1+\alpha)\int_{\R^2}u_n^2{\rm d}x+2\alpha B(u_n)\\
&+\frac{p-4-2\alpha}{p}\lambda\int_{\R^2}|u_n|^p{\rm d}x
+\frac{q-4-2\alpha}{q}\beta\int_{\R^2}|u_n|^q{\rm d}x.
\endaligned
\end{equation}
which, together with (\ref{eqn:fun0}), implies
that $\{u_n\}$ is a bounded sequence in $H_{rad}^1(\R^2)$. In view of Lemma \ref{Lem:xiaozhi}, we have
$$
\aligned
\lim\limits_{n\rightarrow\infty}I_\lambda(u_n)&=\lim\limits_{n\rightarrow\infty}\left(I_{\gamma_n,\beta_n}(u_n)
-\frac{\gamma_n}{4(1+\alpha)}\|u_n\|_4^{4(1+\alpha)}
+\frac{\beta_n}{q}\int_{\R^2}|u_n|^q{\rm d}x\right)\\
&=\lim\limits_{n\rightarrow\infty}c_{\gamma_n,\beta_n}= c^*>\frac{\epsilon^2}{4}.
\endaligned
$$
Moreover, for any $\psi\in C_{0}^\infty(\R^2)\cap H_{rad}^1(\R^2)$, we have
$$
\aligned
&\lim\limits_{n\rightarrow\infty}\lan I_\lambda'(u_n),\psi\ran\\
&=\lim\limits_{n\rightarrow\infty}
\left(\lan I_{\gamma_n,\beta_n}'(u_n),\psi\ran-\gamma_n\|u_n\|_4^{\alpha}\int_{\R^2}|u_n|^2u_n\psi {\rm d}x
+\beta_n\int_{\R^2}|u_n|^{q-2}u_n\psi {\rm d}x\right)=0.
\endaligned
$$
That is to say, $\{u_n\}$ is a bounded Palais-Smale sequence for $I_\lambda$ at level $c^*$.
Thus, there exists $u^*\in H_{rad}^1(\R^2)$ such
that $u_n\rightharpoonup u^*$ weakly in $H_{rad}^1(\R^2)$
and $u_n\rightarrow u^*$ strongly in $L^q(\R^3)$ for $q\in(2,+\infty)$.
Similarly to the arguments of Lemma \ref{Lem:PS}, we see that $I_\lambda'(u^*)=0$
and $u_n\rightarrow u^*$ strongly in $H_{rad}^1(\R^2)$ as $n\rightarrow0$.
Thus, the fact that $u_n\in H_{rad}^1(\R^2)\setminus(P_\epsilon^+\cup P_\epsilon^-)$
yields $u^*\in H_{rad}^1(\R^2)\setminus(P_\epsilon^+\cup P_\epsilon^-)$ and then
$u^*$ is a sign-changing solution of equation (\ref{eqn:cs3}).\\
\textbf{Step 3.} Define
$$
\bar{c}:=\inf\limits_{u\in\Theta}I_\lambda(u),\quad\Theta:=\{u\in H_{rad}^1(\R^2)\setminus\{0\},\,I_\lambda'(u)=0,\,u^{\pm}\not\equiv0\}.
$$
Based on above, we see $\Theta\not=\emptyset$ and
$\bar{c}\leq c^*$, where $c^*$ is given in the Step 2.
 By the definition of $\bar{c}$, there exists
$\{u_n\}\subset H_{rad}^1(\R^2)$ such that $I_\lambda(u_n)\rightarrow \bar{c}$ and $I_\lambda'(u_n)=0$. Using the earlier arguments,
we can deduce that $\{u_n\}$ is a bounded sequence in $H_{rad}^1(\R^2)$.
Arguing as in Lemma \ref{Lem:PS}, there exists a nontrivial
$u\in H_{rad}^1(\R^2)$ such that $I_\lambda(u)=\bar{c}$ and $I_\lambda'(u)=0$.
On the other hand, we deduce from $\langle I_\lambda'(u_n),u_n^\pm\rangle=0$ that for any $\varepsilon>0$,
there exists $C_\varepsilon>0$ such that
$$
C\|u_n^\pm\|_p^2\leq\|u_n^\pm\|^2\leq\lambda\int_{\R^2}|u_n|^{p-2}u_nu_n^\pm {\rm d}x
=\lambda\int_{\R^2}|u_n^\pm|^{p}{\rm d}x
$$
which, together with the boundedness of $\{u_n\}$ in $H_{rad}^1(\R^2)$ implies that $\|u_n^\pm\|_p\geq C$. Hence,
$\|u^\pm\|_p\geq C$, and then we see that
$u$ is a ground state sign-changing solution of equation (\ref{eqn:cs3}). The proof is complete.
\qed


\s{Asymptotic behavior of nodal solutions as $\lambda\rightarrow+\infty$}
\renewcommand{\theequation}{4.\arabic{equation}}

In this section we are attempt to get energy double behavior of the sign-changing
solution obtained in Theorem \ref{Thm:existence} as $\lambda$ tends to infinity.\\
We state the following modified equation:
\begin{equation}\label{eqn:modified}
\left\{
  \begin{array}{ll}
    -\triangle {u}+\omega u+\bar{\lambda}\mathcal{B}(u)u =|u|^{p-1}u& \mbox{in}\,\,\R^2, \\
    u(x)=u(|x|)\,  \in\, H^1(\R^2),
  \end{array}
\right.
\end{equation}
where $\bar{\lambda}=\lambda^{\frac{-4}{p-2}}$ and $\mathcal{B}(u)$ is given in Section 3.
The energy functional $J_{\bar{\lambda}}$ corresponding to equation (\ref{eqn:modified})
is denoted as
$$
J_{\bar{\lambda}}(u)=\frac{1}{2}\int_{\R^2}(|\nabla u|^2+\omega u^2){\rm d}x+\frac{\bar{\lambda}}{2}\int_{\R^2}\frac{u^2}{|x|^2}
h^2(|x|){\rm d}x-\frac{1}{p}\int_{\R^2}|u|^{p}{\rm d}x.
$$
It is easy to check that for any $u\in H_{rad}^1(\R^2)$, the identity
$J_{\bar{\lambda}}(\bar{u})=\lambda^{\frac{2}{p-2}}I_\lambda(u)$ holds,
where $\bar{u}=\lambda^{\frac{1}{p-2}}u$.
In particular, $u$ solves equation (\ref{eqn:cs3}) if and only if $\bar{u}$ satisfies
the modified equation (\ref{eqn:modified}).

\begin{lemma}\label{Lem:modfi}
$w\in H_{rad}^1(\R^2)$ is a ground state sign-changing solution of equation (\ref{eqn:cs3}) if and only if $\bar{w}$ is
a ground state sign-changing solution of equation (\ref{eqn:modified}), where $\bar{w}=\lambda^{\frac{1}{p-2}}w$.
\end{lemma}
\Proof
It is not hard to check that $J_{\bar{\lambda}}(\bar{w})=\lambda^{\frac{2}{p-2}}I_\lambda(w)$.\\
\emph{ Sufficiency:}
Assume on the contrary that $\bar{w}$ is not a ground state sign-changing
solution of equation (\ref{eqn:modified}),
then there exists $\bar{v}\in H_{rad}^1(\R^2)$ is a sign-changing solution
of equation (\ref{eqn:modified}) with
\begin{equation}\label{eqn:bijiao1}
J_{\bar{\lambda}}(\bar{v})<J_{\bar{\lambda}}(\bar{w}).
\end{equation}
It follows from (\ref{eqn:modified}) that there exists sign-changing solution
$v\in H_{rad}^1(\R^2)$ of equation (\ref{eqn:cs3}) such that
\begin{equation}\label{eqn:bijiao2}
J_{\bar{\lambda}}(\bar{v})=\lambda^{\frac{2}{p-2}}I_\lambda(v)\geq
\lambda^{\frac{2}{p-2}}I_\lambda(w)=J_{\bar{\lambda}}(\bar{w}),
\end{equation}
which contradicts (\ref{eqn:bijiao1}). \\
\emph{Necessity:} If $\bar{w}$ is
a ground state sign-changing solution of equation (\ref{eqn:modified}),
then we are attempt to obtain a contradiction by assuming that
$w:=\lambda^{\frac{-1}{p-2}}\bar{w}\in H_{rad}^1(\R^2)$
is not a ground state sign-changing
solution of equation (\ref{eqn:cs3}).
Clearly, $w$ solves equation (\ref{eqn:cs3}). Thus, there exists
$v\in H^1_{rad}(\R^2)$ is a sign-changing solution
of (\ref{eqn:cs3}) with
\begin{equation}\label{eqn:bijiao3}
I_{\lambda}(v)<I_{\lambda}(w).
\end{equation}
In view of equations (\ref{eqn:modified}) and (\ref{eqn:cs3}) that there exists sign-changing solution
$\bar{v}\in H_{rad}^1(\R^2)$ of equation (\ref{eqn:modified}) such that
$$
I_{\lambda}(v)=\lambda^{\frac{-2}{p-2}}J_{\bar{\lambda}}(\bar{v})\geq \lambda^{\frac{-2}{p-2}}J_{\bar{\lambda}}(\bar{w})=I_{\lambda}(w),
$$
which is a contradiction.
The proof is complete.
\qed
\vskip0.1in
Now we prove Theorem \ref{Thm:ground}. To emphasize the dependence on $\lambda$, we use
$w_\lambda\in H_{rad}^1(\R^2)$ to denote a ground state sign-changing solution of (\ref{eqn:cs3}) obtained in
Theorem \ref{Thm:existence}. Then by Lemma {\ref{Lem:modfi}}, we know that $\bar{w}_{\bar{\lambda}}$
is also a ground state sign-changing solution of (\ref{eqn:modified}).

\textbf{Proof of Theorem \ref{Thm:ground}}.

Define
\begin{equation}\label{eqn:mb}
m_{\bar{\lambda}}:=\inf\limits_{u\in\Theta}J_{\bar{\lambda}}(u),
\quad\Theta:=\{u\in H_{rad}^1(\R^2)\setminus\{0\},\,J_{\bar{\lambda}}'(u)=0,\,u^{\pm}\not\equiv0\}.
\end{equation}
In view of Theorem \ref{Thm:existence} and Lemma \ref{Lem:modfi},
we know $\Theta\not=\emptyset$. Moreover, if we replace $I_\lambda$ by $J_{\bar{\lambda}}$,
then all the above calculations in Sect. 3 can be repeated word by word.
So we can also get
a sign-changing solution to equation (\ref{eqn:modified}).
Furthermore, we can deduce from (\ref{eqn:fun0}) that
\begin{equation}\label{eqn:fun0--}
m_{\bar{\lambda}}\leq C_R:=\sup\limits_{u\in\psi_0(\triangle)}\bar{J}_{1}(u)<\infty,
\end{equation}
where $\psi_0(\triangle)$ has been defined in the proof of Theorem \ref{Thm:existence} and
$C_R$ is independent of $\bar{\lambda}\in(0,1]$ and
the functional $\bar{J}_{1}: H_{rad}^1(\R^2)\rightarrow\R$ is defined as
$$
\bar{J}_{1}:=\frac{1}{2}\|u\|^2+\frac{1}{2}\int_{\R^2}\frac{u^2}{|x|^2}
h^2(|x|){\rm d}x+\frac{1}{4(1+\alpha)}\|u\|_4^{4(1+\alpha)}-\frac{1}{p}\int_{\R^2}|u|^{p}{\rm d}x.
$$
We now show that, for $\bar{\lambda}\rightarrow0^+$, $\{\bar{w}_{\bar{\lambda}}\}_{\bar{\lambda}\in(0,1)}$
is bounded in $H_{rad}^1(\R^2)$.
 Indeed,
\begin{equation}\label{eqn:fun1--}
m_{\bar{\lambda}}=\frac{1}{2}\int_{\R^2}(|\nabla \bar{w}_{\bar{\lambda}}|^2{\rm d}x
+\omega \bar{w}_{\bar{\lambda}}^2){\rm d}x
+\frac{\bar{\lambda}}{2}\int_{\R^2}\frac{\bar{h}^2(|x|)}{|x|^2} \bar{w}_{\bar{\lambda}}^2{\rm d}x
-\frac{1}{p}\int_{\R^2}|\bar{w}_{\bar{\lambda}}|^p{\rm d}x,
\end{equation}
where $\bar{h}(|x|)=\int_0^{|x|}\frac{\tau}{2} \bar{w}_{\bar{\lambda}}^2(\tau){\rm d}\tau$, and
 \begin{equation}\label{eqn:fun2--}
0=\int_{\R^2}|\nabla \bar{w}_{\bar{\lambda}}|^2{\rm d}x+\int_{\R^2}\omega \bar{w}_{\bar{\lambda}}^2{\rm d}x
+3\bar{\lambda}\int_{\R^2}\frac{\bar{h}^2(|x|)}{|x|^2} \bar{w}_{\bar{\lambda}}^2{\rm d}x
-\int_{\R^2}|\bar{w}_{\bar{\lambda}}|^p{\rm d}x.
\end{equation}
Moreover, we have the following Pohozaev identity
\begin{equation}\label{eqn:fun3--}
0=\omega\int_{\R^2}\bar{w}_{\bar{\lambda}}^2{\rm d}x
+2\bar{\lambda}\int_{\R^2}\frac{\bar{h}^2(|x|)}{|x|^2} \bar{w}_{\bar{\lambda}}^2{\rm d}x
-\frac{2}{p}\int_{\R^2}|\bar{w}_{\bar{\lambda}}|^p{\rm d}x.
\end{equation}
Multiplying (\ref{eqn:fun1--}) by $2$, (\ref{eqn:fun2--}) by $-1$ and (\ref{eqn:fun3--}) by $1$ and adding them up,
we get
\begin{equation}\label{eqn:fun4--}
2m_{\bar{\lambda}}=\int_{\R^2}\omega \bar{w}_{\bar{\lambda}}^2{\rm d}x
+\frac{p-4}{p}\int_{\R^2}|\bar{w}_{\bar{\lambda}}|^p{\rm d}x.
\end{equation}
Based on conclusion (i) of Lemma \ref{Lem:hs} and interpolation inequality, we have
$$
B(\bar{w}_{\bar{\lambda}})=\int_{\R^2}\frac{\bar{h}^2(|x|)}{|x|^2}\bar{w}_{\bar{\lambda}}^2{\rm d}x
\leq C\|\bar{w}_{\bar{\lambda}}\|_4^4\|\bar{w}_{\bar{\lambda}}\|_2^2
\leq\|\bar{w}_{\bar{\lambda}}\|_2^{\frac{2(p-4)}{p-2}}\|\bar{w}_{\bar{\lambda}}\|_p^{\frac{2p}{p-2}}
\|\bar{w}_{\bar{\lambda}}\|_2^2,
$$
which, together with (\ref{eqn:fun0--}), (\ref{eqn:fun1--}) and (\ref{eqn:fun4--}), implies
that $\{\bar{w}_{\bar{\lambda}}\}_{\bar{\lambda}\in(0,1)}$ is bounded in $H_{rad}^1(\R^2)$.
Up to a subsequence, we assume that $\bar{w}_{\bar{\lambda}}
\rightarrow w_0$ weakly in $H_{rad}^1(\R^2)$ as $\bar{\lambda}\rightarrow0^+$.
Note that $\bar{w}_{\bar{\lambda}}$ is a ground state sign-changing solution
of equation (\ref{eqn:modified}),
then, similarly to the arguments of Lemma \ref{Lem:PS}, we deduce that
$\bar{w}_{\bar{\lambda}}
\rightarrow \bar{w}_0$ strongly in $H_{rad}^1(\R^2)$ and $\bar{w}_0$ is a sign-changing solution
of equation (\ref{eqn:ellptic}). Thus,
we have
\begin{equation}\label{eqn:fun5--0}
m_{\bar{\lambda}}=J_{\bar{\lambda}}(\bar{w}_{\bar{\lambda}})
=I(\bar{w}_{0})+o(\bar{\lambda})\geq m_0+o(\bar{\lambda}),
\end{equation}
where $m_0$ is denoted as a ground state sign-changing energy of equation (\ref{eqn:ellptic}).
By the energy doubling property of sign-changing solution of equation (\ref{eqn:ellptic}), we get
$m_0>2c_0$, where $c_0$ is defined in (\ref{eqn:p-ground}).  In view of (\ref{eqn:fun5--0}), we have
\begin{equation}\label{eqn:fun5--}
m_{\bar{\lambda}}=I(\bar{w}_{0})+o(\bar{\lambda})=m_0+o(\bar{\lambda})>2c_0
\end{equation}
for $\bar{\lambda}$ small enough.

Let $u_\lambda\in\mathcal{N}_\lambda$ be a positive ground
state solution of equation (\ref{eqn:cs3}) (see \cite{Byeon12}),
then $I_\lambda(u_\lambda)=c_\lambda$, where
$\mathcal{N}_\lambda$ and $c_\lambda$ have been defined in (\ref{eqn:p-ground+})
and (\ref{eqn:p-ground1}), respectively. Similar to that argument in Lemma \ref{Lem:modfi},
we can prove that there exists a positive ground
state solution $\bar{u}_{\bar{\lambda}}:=\lambda^{\frac{1}{p-2}}u_\lambda$
of equation (\ref{eqn:modified}) with $J_{\bar{\lambda}}(\bar{u}_{\bar{\lambda}})
=c_{\bar{\lambda}}$. Obviously, $c_{\bar{\lambda}}=\lambda^{\frac{2}{p-2}}c_\lambda$
Furthermore, taking $\bar{\lambda}\rightarrow0^+$, we
also prove that there exists $\bar{u}_0\in H_{rad}^1(\R^2)$ such that $\bar{u}_{\bar{\lambda}}
\rightarrow \bar{u}_0$ strongly in $H_{rad}^1(\R^2)$ and $\bar{u}_0$ is a
positive solution of equation (\ref{eqn:ellptic}).
Hence, by the uniqueness of positive solution of (\ref{eqn:ellptic}), we have
\begin{equation}\label{eqn:fun6--}
c_{\bar{\lambda}}=J_{\bar{\lambda}}(\bar{u}_{\bar{\lambda}})=I(\bar{u}_0)+o(\bar{\lambda})=c_0+o(\bar{\lambda}),
\end{equation}
where $c_0$ is given in (\ref{eqn:p-ground}). Combining (\ref{eqn:fun5--}) with (\ref{eqn:fun6--}), we know that there exists $\bar{\lambda}^*>0$
such that $m_{\bar{\lambda}}>2c_{\bar{\lambda}}$ for any $\bar{\lambda}\in(0,\bar{\lambda}^*)$.
It is easy to check that $m_{\lambda}$ is strictly larger than twice
that of the ground state energy $c_\lambda$ for $\lambda\in(\lambda^*,+\infty)$, where $\lambda^*=(\bar{\lambda}^*)^\frac{2-p}{4}$.
The proof is complete.
\qed


\s{Multiplicity}
\renewcommand{\theequation}{5.\arabic{equation}}
In this section, we prove the existence of infinitely many sign-changing
 solutions to equation (\ref{eqn:cs3}) when $\lambda$ is small enough.
\subsection{The perturbed problem}
We here employ a perturbed approach which is introduced in \cite{Liuj15,Liuz16}
to overcome the difficulty getting the boundedness of the Palais-Smale,
due to lack of the well-known {\it Ambrosetti-Rabinowtiz} condition.
For any fixed $\beta\in(0,1]$ and $q\in(6,8)$, we consider the modified problem
\begin{equation}\label{eqn:4-peturbed}
-\triangle u+ \omega u+\mathcal{B}(u)u\\
=\lambda|u|^{p-2}u+\beta|u|^{q-2}u,\quad x\in\R^2
\end{equation}
and its associated functional is given as below
$$
I_\beta(u)=I_\lambda(u)-\frac{\beta}{q}\int_{\R^3}|u|^q{\rm d}x,
$$
where $\mathcal{B}(u)$ is given in Section 3.
It is easy to show that $I_\beta\in C^1(H_{rad}^1(\R^2),\R)$ and
$$
 \langle I_\beta'(u),v\rangle=  \langle I_\lambda'(u),v \rangle
 -\beta\int_{\R^3}|u|^{r-2}uv{\rm d}x,\,\,u,v\in H_{rad}^1(\R^2).
$$
Similarly to (\ref{eqn:in1}), for each $u\in H_{rad}^1(\R^2)$, the following equation
\begin{equation} \label{eqn:4-in1}
-\triangle v+[\omega+\mathcal{B}(u)] v
=\lambda|u|^{p-2}u+\beta|u|^{q-2}u,\quad x\in\R^2
\end{equation}
has a unique weak solution $v\in H_{rad}^1(\R^2)$. In order to construct the descending flow
for functional $I_\beta$,
we also introduce an auxiliary operator $T_\beta: u\in H_{rad}^1(\R^2)\mapsto v\in H_{rad}^1(\R^2)$, where $v=T_\beta(u)$ is the unique weak solution of problem (\ref{eqn:in1}).
Clearly, the fact that $u$ is a solution of problem (\ref{eqn:4-in1}) is equivalent to that $u$ is a fixed point of $T_\beta$.
As in Sect.2, one can obtain that the operator $T_\beta: H_{rad}^1(\R^2)\rightarrow H_{rad}^1(\R^2)$ is well defined and is continuous.
In the following, if the proof of a result is similar to its counterpart in Section 2,
it will not be written out.

\begin{lemma}\label{Lem:4-T2-4}
\begin{itemize}
\item[\rm (1) ] $\lan I'_{\beta}(u), u-T_{\beta}(u)\ran\geq \|u-T_{\beta}(u)\|^2$ for all $u\in H_{rad}^1(\R^2)$;
\item[\rm (2) ] $\|I'_{\beta}(u)\|\leq \|u-T_{\beta}(u)\|(1+C\|u\|^2)$ for some $C>0$ and all
$u\in H_{rad}^1(\R^2)$.
\end{itemize}
\end{lemma}

\begin{lemma}\label{Lem:4-T3-}
For any fixed $\beta\in(0,1]$, $c<d$ and $\alpha>0$, there exists $\delta>0$ (which depends on $\beta$) and $\lambda_*>0$ such that $\|u-T_{\beta}(u)\|\geq \delta$
if $u\in H_{rad}^1(\R^2)$, $I_{\beta}(u)\in [c,d]$
and $\|I_{\beta}'(u)\|\geq\alpha$ and $\lambda\in(0,\lambda_*)$.
\end{lemma}
\Proof
Let $\mu\in(6,q)$. For any $u\in H_{rad}^1(\R^2)$, we have
\begin{equation}\label{eqn:4-T0}
\aligned
&I_{\beta}(u)-\frac{1}{\mu}\langle u,u-T_{\beta}(u)\rangle\\
&=\frac{\mu-2}{2\mu}\|u\|^2
+\frac{1}{\mu}\int_{\R^2}\mathcal{B}(u)u(u-T_{\beta}(u)){\rm d}x
+\frac{\mu-2}{\mu}B(u)\\
&\,\,\,\,\,\,\,+\frac{(p-\mu)\lambda}{\mu p}\int_{\R^2}|u|^{p}{\rm d}x-\frac{1}{\mu}\int_{\R^2}u^2\int_{|x|}^{+\infty}\frac{h(s)}{s}u^2(s){\rm d}s{\rm d}x
+\frac{(q-\mu)\beta}{\mu q}\int_{\R^2}|u|^{q}{\rm d}x.
\endaligned
\end{equation}
Then, by the conclusion (iv) of Lemma \ref{Lem:bu}, we have
\begin{equation}\label{eqn:4-T1}
\aligned
&\frac{\mu-2}{2\mu}\|u\|^2+\frac{\mu-6}{\mu}B(u)+
\frac{(q-\mu)\beta}{\mu q}\beta\|u\|_q^q+\frac{(p-\mu)\lambda}{\mu p}\|u\|_{p}^{p}\\
&\le C\left(|I_{\gamma,\beta}(u)|+\|u\|\|u- T_{\gamma,\beta}(u)\|
+\int_{\R^2}\mathcal{B}(u)|u(u-T_{\gamma,\beta}(u))|{\rm d}x\right).
\endaligned
\end{equation}
By Lemma \ref{Lem:inequality}, Young's inequality and (\ref{eqn:4-T1}) we have
\begin{equation}\label{eqn:4-T2}
\aligned
&\frac{\mu-6}{8\mu}\|u\|_4^4+\frac{\mu-2}{4\mu}\|u\|^2+\frac{\mu-6}{2\mu}B(u)
+\frac{(q-\mu)\beta}{\mu q}\|u\|_q^q+\frac{(p-\mu)\lambda}{\mu p}\|u\|_{p}^{p}\\
&\le C\left(|I_{\beta}(u)|+\|u\|\|u- T_{\beta}(u)\|
+\int_{\R^2}\mathcal{B}(u)|u(u-T_{\beta}(u))|{\rm d}x\right).
\endaligned
\end{equation}
Set $r=q$ in the conclusion (iii) of Lemma \ref{Lem:hs},
then by Lemma \ref{Lem:hs} and H\"{o}lder's inequality, for any $\xi>0$, there exists $C_\xi>0$ such that
\begin{align}\label{eqn:4-T3}
&\int_{\R^2}|u(u-T_{\beta}(u))|\int_{|x|}^{+\infty}\frac{h(s)}{s}u^2(s){\rm d}s{\rm d}x\nonumber\\
&=2\pi\int_{0}^\infty|u(\tau)(u(r)-T_{\beta}(u(\tau)))\tau|\int_{\tau}^{+\infty}\frac{h(s)}{s}u^2(s){\rm d}s{\rm d}\tau\nonumber\\
&=2\pi\int_{0}^{+\infty}\frac{h(s)}{s}u^2(s){\rm d}s\int_{0}^s|u(\tau)(u(r)-T_{\beta}(u(\tau)))\tau|{\rm d}\tau\nonumber\\
&\leq2\pi\int_{0}^{+\infty}\frac{h(s)}{s}u^2(s){\rm d}s
\left[\xi\int_{0}^su^2(\tau)\tau{\rm d}\tau+C_{\xi}\int_{0}^s|(u(r)-T_{\beta}(u(\tau)))^2\tau{\rm d}\tau\right]\nonumber\\
&\leq2\pi\xi\int_{0}^{+\infty}\frac{h(s)}{s}u^2(s){\rm d}s\int_{0}^su^2(\tau)\tau{\rm d}\tau
+2\pi C_{\xi}\int_{0}^{+\infty}\frac{h(s)}{s}u^2(s){\rm d}s\int_{0}^s(u(r)-T_{\beta}(u(\tau)))^2\tau{\rm d}\tau\nonumber\\
&\leq4\pi\xi\int_{0}^{+\infty}\frac{h^2(s)}{s}u^2(s){\rm d}s+C_{\xi}\|u-T_{\beta}(u)\|_2^2\int_{0}^{+\infty}\frac{h(s)}{s}u^2(s){\rm d}s\nonumber\\
&\leq4\xi B(u)+C_{\xi}\|u-T_{\beta}(u)\|_2^2\left[\frac{2}{\pi}B(u)+\frac{1}{2\pi}\|u\|^2_2
+C(\|u\|_{4}^4+\|u\|_q^{\frac{q}{2}})\right].
\end{align}
Assume on the contrary that there exists $\{u_n\}\subset H_{rad}^1(\R^2)$
with $I_{\beta}(u_n)\in[c,d]$ and
$\|I'_{\beta}(u_n)\|\geq \alpha$ such that
$$
\|u_n-T_{\beta}(u_n)\|\rightarrow0,\quad\text{as}\,\,n\rightarrow\infty.
$$
Combining (\ref{eqn:T3-0}), (\ref{eqn:4-T2}) and (\ref{eqn:4-T3}), we have for large $n$
\begin{equation}\label{eqn:4-T4}
\aligned
&\frac{\mu-6}{8\mu}\|u_n\|_4^4+\frac{\mu-2}{4\mu}\|u_n\|^2+\frac{\mu-6}{2\mu}B(u_n)
+\frac{(q-\mu)\beta}{\mu q}\|u_n\|_q^q+\frac{(p-\mu)\lambda}{\mu p}\|u_n\|_{p}^{p}\\
&\le C[|I_{\beta}(u_n)|+\|u_n\|\|u_n- T_{\beta}(u_n)\|
+\|u_n\|_q^{\frac{q}{2}}\|u_n- T_{\beta}(u_n)\|^2].
\endaligned
\end{equation}
Now we claim that sequence $\{u_n\}$ is bounded in $H_{rad}^1(\R^2)$.
Otherwise, we assume $\|u_n\|\rightarrow\infty$.
It follows from (\ref{eqn:4-T4}) that
\begin{equation}\label{eqn:4-T5}
C_1\|u_n\|_4^4+C_2\|u_n\|^2
+C_3\|u_n\|_q^q-C_4\lambda\|u_n\|_{p}^{p}\\
\le C.
\end{equation}
for large $n$, where $C_i, i=1,2,3,4$ are some positive constants.
Obverse that, there exists $\lambda^*>0$ such that for any $\lambda\in(0,\lambda^*)$,
we have
$$
C_1|t|^4+C_3|t|^q-C_4\lambda |t|^{p}\ge0,\quad t\in\R,
$$
due to $4<p<q$. Applying this with $t=|u_n|$, then we see from (\ref{eqn:4-T5})
that sequence $\{u_n\}$ is bounded in $H_{rad}^1(\R^2)$ for any fixed $\beta\in(0,1]$.
The claim combined with Lemma \ref{Lem:4-T2-4} implies $\|I'_{\beta}(u_n)\|\rightarrow0$ as $n\rightarrow\infty$, which is a contradiction.
\qed

\subsection{Invariant Subsets of Descending Flows}

\begin{lemma}\label{Lem:4-C1-4}
 There exists $\epsilon_0>0$ such that for $\epsilon\in(0,\epsilon_0)$,
\item[\rm (1) ] $T_{\beta}(\partial P_\epsilon^-)\subset P_\epsilon^-$ and every nontrivial solution $u\in P_\epsilon^-$ is negative,
\item[\rm (2) ] $T_{\beta}(\partial P_\epsilon^+)\subset P_\epsilon^+$ and every nontrivial solution $u\in P_\epsilon^+$ is positive.
\end{lemma}

\begin{lemma}\label{Lem:4-local4}
There exists a locally Lipschitz continuous operator $B_{\beta}:\,H_{rad}^1(\R^2)\setminus K_{\beta}\rightarrow H_{rad}^1(\R^2)$ such that
\begin{itemize}
\item[(i)] $\langle T'_{\beta}(u),u-B_{\beta}(u)\rangle\geq\frac{1}{2}\|u-T_{\beta}(u)\|^2$;
\item[ (ii)]$\frac{1}{2}\|u-B_{\beta}(u)\|^2\leq\|u-T_{\beta}(u)\|^2\leq2\|u-B_{\beta}(u)\|^2$;
\item[ (iii)]$T_{\beta}(\partial P_\epsilon^\pm)\subset P_\epsilon^\pm$, $\forall \epsilon\in(0,\epsilon_0)$;
\item[ (iv)] $B_{\beta}$ is odd,
\end{itemize}
where $K_{\beta}:=\{u\in H_{rad}^1(\R^2):\,I_\beta'(u)=0 \}$ and $\epsilon_0$ has been given in Lemma \ref{Lem:4-C1-4}.
\end{lemma}
Moreover, we can prove the functional $I_{\beta}$ satisfies the (PS)-condition
with the aid of Lemma \ref{Lem:4-T3-}.

In order to obtain infinitely many sign-changing solutions, we will make use of an abstract critical point
developed by Liu et al \cite{Liuj15}, which we
recall below.
The notations from Section 2 are still valid. Assume $G : E\rightarrow E$ to be an
isometric involution, that is, $G^2=id$ and $d(Gx;Gy) = d(x;y)$ for $x,y\in E$. We
assume $J$ is $G$-invariant on $H_{RAD}^1(\R^2)$ in the sense that $J(Gx) = J(x)$ for any $x\in E$.
We also assume $Q = GP$. A subset $F\subset E$ is said to be symmetric if $Gx \in F$
for any $x \in F$. The genus of a closed symmetric subset $F$ of $E\setminus\{0\}$ is denoted
by $\gamma(F)$.

\begin{definition}{\rm \cite{Liuj15}}\label{Def:2}
$P$ is called a $G$-admissible invariant set with respect to
$J$ at level $c$, if the following deformation property holds: there exist $\epsilon_0>0$ and
a symmetric open neighborhood $N$ of $K_c\setminus{W}$ with $\gamma(\bar{N})< \infty$, such that for
$\epsilon\in(0,\epsilon_0)$ there exists $\eta\in C(E,E)$ satisfying
\begin{itemize}
\item[\rm (1) ]
$\eta(\bar{P})\subset\bar{P}$, $\eta(\bar{Q})\subset\bar{Q}$;
\item[\rm (2) ]$\eta\circ G=G\circ\eta$;
\item[\rm (2) ]$\eta|_{J^{c-2\epsilon}}=id$;
\item[\rm (3) ] $\eta(J^{c+\epsilon}\setminus{(N\cup W})\subset J^{c-\epsilon}$.
\end{itemize}
\end{definition}

\begin{theorem}{\rm \cite{Liuj15}}\label{Thm:xjl2}
 Assume that $P$ is a $G$-admissible invariant set with respect
to $J$ at any level $c\geq c^*:= \inf_{u\in\Sigma} J(u)$ and for any $n\in \mathbb{N}$, there exists a
continuous map $\psi_k: B_k := \{x \in \R^k : |x|\leq1\}\rightarrow E$ satisfying
\begin{itemize}
\item[\rm (1) ] $\psi_k(0)\subset M:=P\cap Q$ and $\psi_k(-t)=G \psi_k(t)$ for $t\in B_k$,
\item[\rm (2) ]$\psi_k(\partial B_k)\cap M=\emptyset$,
\item[\rm (3) ] $\sup\limits_{u\in Fix_{G}\cup\psi_k(\partial B_k)}J(u)<c_*$,
where $Fix_{G}:=\{u\in E;Gu=u\}$.
\end{itemize}
For $j\in N$, define
$$
c_j=\inf\limits_{B\in\Gamma_j}\sup\limits_{u\in B\setminus{W}}J(u),
$$
where
$$
\aligned
\Gamma_j:=\{&B\big{|}B=\psi(B_k\setminus{Y})for\, some\,\psi\in G_k, k\geq j,\\
&and\, open Y\subset B_k \,such\, that
-Y=Y\,and\, \gamma(\bar{Y})\leq k-j\}
\endaligned
$$
and
$$
G_k:=\{\psi|\psi\in C(B_k,E), \psi(-t)=G\psi(t)\,for\,t\in B_k,\psi(0)\in M\,
and\,\psi|_{\partial B_k}=\psi_k|_{\partial B_k}
\}.
$$
Then for $j\geq2$, $c_j\geq c_*$ and $K_{c_j}\setminus{W}\not=\emptyset$ and $c_j\rightarrow\infty$ as $j\rightarrow\infty$.
\end{theorem}

In order to apply Theorem \ref{Thm:xjl2}, we set $E=H_{rad}^1(\R^2)$, $G = -id$, $J = I_\beta$ and $P = P_\epsilon^+$. Then
$M = P_\epsilon^+
\cap P_\epsilon^-$, $\Sigma=\partial P_\epsilon^+
\cap \partial P_\epsilon^-$, and $W = P_\epsilon^+\cup P_\epsilon^-$.
 Now, we show that $P_\epsilon^+$
is a G-admissible invariant set for the functional $I_\beta$ at any level $c$. Since $K_c$ is
compact, there exists a symmetric open neighborhood $N$ of $K_c \setminus{W}$ such that
$\gamma(\bar{N})<\infty$. Similar to Lemma 3.9 in \cite{Liuz16}, we have

\begin{lemma}\label{Lem:deform1}
There exists $\epsilon_0 > 0$ such that for $0 <\epsilon<\epsilon'<\epsilon_0$, there exists a
continuous map $\eta : [0,1]\times H_{rad}^1(\R^2) \rightarrow H_{rad}^1(\R^2)$ satisfying
\begin{itemize}
\item[\rm (1) ] $\eta(0,u)=u$ for $u\in H_{rad}^1(\R^2)$.
\item[\rm (2) ] $\eta(t,u)=u$ for $t\in[0,1]$, $u\not\in I_\beta^{-1}[c-\epsilon',c+\epsilon']$.
\item[\rm (3) ] $\eta(t,-u)=-\eta(t,u)$ for $(t,u)\in[0,1]\times H_{rad}^1(\R^2)$.
\item[\rm (4) ]$\eta(1,I_\beta^{c+\epsilon}\setminus{(N\cup W)})\subset I_\beta^{c-\epsilon}$.
\item[\rm (5) ] $\eta(t,\overline{P_\epsilon^+})\subset\overline{P_\epsilon^+}$,
$\eta(t,\overline{P_\epsilon^-})
\subset\overline{P_\epsilon^-}$
for $t\in[0,1]$.
\end{itemize}
\end{lemma}
\subsection{Proof of Theorem \ref{Thm:many}(Multiplicity)}
 We use two steps to give the proof.\\
\textbf{Step1.} It follows from Lemma \ref{Lem:deform1}
that $P_\epsilon^+$ is a G-admissible invariant set
for the functional $I_{\beta}$ for $\beta\in(0,1]$ at any level $c$.
We construct $\psi_k$ satisfying the hypotheses of Theorem \ref{Thm:xjl2}.
 For any fixed $k\in\mathbb{N}$, we choose
$\{v_i\}_{i=1}^k\subset \{C_{0}^\infty(\R^2)\cap H_{rad}^1(\R^2)\}\setminus\{0\}$ such that $\text{supp}(v_i)\cap\text{supp}(v_j)$
for $i\not=j$. Define
$\psi_k\in C(B_k,H_{rad}^1(\R^2))$ as
$$
\psi_k(t)(\cdot)=R_k\Sigma_{i=1}^kt_iv_i(R_k\cdot),\quad t=(t_1,t_2,...,t_k)\in B_k.
$$
Observe that
$$
\rho_k=\min\{\|t_1v_1+t_2v_2+\cdot\cdot\cdot+t_kv_k\|^2_2:\,\,\sqrt{\Sigma_{i=1}^k t_i^2}= 1\}>0,
$$
then $\|u_t\|_2^2\geq\rho_k $ for $u\in\psi_k(\partial B_k)$ and
it follows from Lemma \ref{Lem:qfanshu} that
$\psi_k(\partial B_k)\cap P_\epsilon^+\cap P_\epsilon^-=\emptyset$. Similarly to the proof of
Theorem \ref{Thm:existence},
we also have
$$
\sup\limits_{u\in\psi_k(\partial B_k)}I_{\beta}(u)<0<\inf\limits_{u\in\Sigma} I_{\beta}(u).
$$
Clearly, $\psi_k(0)=0\in P_\epsilon^+\cap P_\epsilon^-$ and $\psi_k(-t)=-\psi_k(t)$ for $t\in B_k$.
For any fixed $\beta\in(0,1)$ and $j\in\{1,2,...,k\}$, we define
$$
c_{\beta}^j=\inf\limits_{B\in\Gamma_j}\sup\limits_{u\in B\setminus{W}}I_{\beta}(u),
$$
where $\Gamma_j$ has been defined in Theorem \ref{Thm:xjl2}.
Based on Theorem \ref{Thm:xjl2}, for any fixed $\beta\in(0,1]$ and $j\geq2$,
\begin{equation}\label{eqn:fun9}
\frac{\epsilon^2}{4}\leq\inf\limits_{u\in\Sigma}I_{\beta}(u)
:=c^*_{\beta}\leq c_{\beta}^j\rightarrow\infty,\quad \text{as}\,\,j\rightarrow\infty
\end{equation}
and
there exists $\{u^j_{\beta}\}\subset H_{rad}^1(\R^2)\setminus{W}$ such that
$I_{\beta}(u^j_{\beta})=c_{\beta}^j$ and $I'_{\beta}(u^j_{\beta})=0$.\\
\textbf{Step2.} Using the similar way as that arguments of Theorem \ref{Thm:existence}(existence part), for any fixed
$j\geq 2$, $\{u^j_{\beta}\}_{\beta\in(0,1]}$ is bounded in $H_{rad}^1(\R^2)$, that is to say,
there exists
$C>0$ independent of $\beta$ such that $\|u^j_{\beta}\|\leq C$. Without loss of generality,
we assume that
$u^j_{\beta}\rightharpoonup u^j_*$ weakly in $H_{rad}^1(\R^2)$ as $\beta\rightarrow0^+$.
By Lemma \ref{Lem:xiaozhi} and Theorem \ref{Thm:xjl2} we define
$$
\frac{\epsilon^2}{4}\leq
\inf\limits_{u\in\Sigma} I_{\beta}(u)\leq c_{\beta}^j\leq c_{R_k}:=\sup\limits_{u\in\psi_k(B_k)}
I(u),
$$
where $c_{R_k}$ is independent of $\beta$. Let $c_{\beta}^j\rightarrow c_*^j$
as $\beta\rightarrow0^+$. Then we can prove that
$u^j_{\beta}\rightarrow u^j_*$ strongly in $H_{rad}^1(\R^2)$ as $\beta\rightarrow0^+$
and $u^k_*\in H_{rad}^1(\R^2)\setminus{W}$
such that
$I'(u^j_*)=0$ and $I(u^j_*)=c_*^j$. We claim that $c_*^j\rightarrow\infty$
as $j\rightarrow\infty$. Indeed,
$c_{\beta}^j$ is decreasing in $\beta$.
Obviously, $c_{\beta}^j\leq c_*^j$ and by (\ref{eqn:fun9}), we have
$c_*^j\rightarrow+\infty$ as $j\rightarrow\infty$.
Therefore, problem (\ref{eqn:cs3}) has infinitely many sign-changing solutions. The proof is complete.
\qed

\end{document}